# An iterative local updating ensemble smoother for estimation and uncertainty assessment of hydrologic model parameters with multimodal distributions

Jiangjiang Zhang[1], Guang Lin[2*], Weixuan Li[3], Laosheng Wu[4], and Lingzao Zeng[1*]

[1] Zhejiang Provincial Key Laboratory of Agricultural Resources and Environment, Institute of Soil and Water Resources and Environmental Science, College of Environmental and Resource Sciences, Zhejiang University, Hangzhou, China,

[2] Department of Mathematics and School of Mechanical Engineering, Purdue University, West Lafayette, Indiana, USA,

[3] Pacific Northwest National Laboratory, Richland, Washington, USA,

[4] Department of Environmental Sciences, University of California, Riverside, California, USA.

* Correspondence to:

    G. Lin, guanglin@purdue.edu

    L. Zeng, lingzao@zju.edu.cn

# Abstract

Ensemble smoother (ES) has been widely used in inverse modeling of hydrologic systems. However, for problems where the distribution of model parameters is multimodal, using ES directly would be problematic. One popular solution is to use a clustering algorithm to identify each mode and update the clusters with ES separately. However, this strategy may not be very efficient when the dimension of parameter space is high or the number of modes is large. Alternatively, we propose in this paper a very simple and efficient algorithm, i.e., the iterative local updating ensemble smoother (ILUES), to explore multimodal distributions of model parameters in nonlinear hydrologic systems. The ILUES algorithm works by updating local ensembles of each sample with ES to explore possible multimodal distributions. To achieve satisfactory data matches in nonlinear problems, we adopt an iterative form of ES to assimilate the measurements multiple times. Numerical cases involving nonlinearity and multimodality are tested to illustrate the performance of the proposed method. It is shown that overall the ILUES algorithm can well quantify the parametric uncertainties of complex hydrologic models, no matter whether the multimodal distribution exists.

# 1. Introduction

Parameter identification is an important aspect in uncertainty quantification of hydrologic systems. However, a direct measurement of model parameters is usually difficult or even impossible in many cases. In this situation, to obtain an estimate of the model parameters, we need to solve an inverse problem with the information provided by some indirect measurements (hereinafter referred to as measurements). Nowadays, Bayesian inversion methods are receiving popularity in hydrologic sciences. In the Bayesian framework, the uncertainties in parameter estimation are represented by the posterior distribution, from which we can obtain any desired statistics [*Stuart*, 2010]. According to Bayes' theorem, the posterior distribution is proportional to the product of the prior distribution times the likelihood. Except for a few simple cases, the analytical form of the posterior distribution is non-existent. In this situation, we have to resort to Monte Carlo simulation methods to sample from the posterior distribution and obtain a numerical approximation accordingly.

One popular method to sample from the posterior distribution is Markov chain Monte Carlo (MCMC), which was first introduced by *Metropolis et al.* [1953] and then extended to more general situations by *Hastings* [1970]. Over the last two decades, many efforts have been devoted to developing efficient MCMC algorithms, including single-chain and multi-chain methods. One of the most popular single-chain MCMC is the delayed rejection adaptive metropolis (DRAM) algorithm developed by *Haario et al.* [2006], which combines the strength of delayed rejection [*Tierney and Mira*, 1999] and adaptive Metropolis [*Haario et al.*, 2001] algorithms. However, when the posterior distribution is multimodal, the performance of the single-chain MCMC would deteriorate [*Vrugt*, 2016]. Through running multiple chains in parallel, MCMC can better explore complex posterior distributions that have multiple modes. One famous example of multi-chain MCMC is the differential evolution adaptive metropolis (DREAM) algorithm [*Vrugt et al.*, 2008; *Vrugt et al.*, 2009b], which is based on the differential evolution Markov chain algorithm [*Braak*, 2006] but uses outlier chain correction and subspace sampling. Due to its efficiency, DREAM has found widespread applications in many different fields [*Vrugt*, 2016]. To sufficiently explore the posterior distribution of model parameters, MCMC usually needs a very large number of model evaluations, especially when the dimension of the parameter space is high. When the system model is CPU-demanding, the computational cost of MCMC simulation would be prohibitive. In this situation, a CPU-efficient surrogate is usually used to replace the original model in MCMC simulation. To eliminate the error introduced by the surrogate, one has to construct an accurate enough surrogate (at least around the posterior distribution [*Zhang et al.*, 2013; *Zhang et al.*, 2016]), or use the original model to correct the surrogate simulation in a two-stage manner [*Efendiev et al.*, 2005; *Laloy et al.*, 2013; *Zeng et al.*, 2012; *Zhang et al.*, 2015]. For the reason of computational cost, it is also difficult to construct an accurate surrogate for a high-dimensional model, except when the nonlinearity of the original model is low enough to allow for a linear approximation [*Li et al.*, 2016b; *Zhang et al.*, 2017].

For parameter estimation in nonlinear problems, a computationally appealing alternative is ensemble Kalman filter (EnKF), which is a Monte Carlo variant of the classical Kalman filter [*Kalman*, 1960]. Since its introduction by *Evensen* [1994], EnKF has been widely used in uncertainty quantification of nonlinear problems in oceanic [*Bertino et al.*, 2003; *Keppenne and Rienecker*, 2003], atmospheric [*Houtekamer and Zhang*, 2016; *Houtekamer and Mitchell*, 2001; *Ott et al.*, 2004], geophysical [*Aanonsen et al.*, 2009; *Gu and Oliver*, 2007] and hydrological [*Chen and Zhang*, 2006; *Moradkhani et al.*, 2005; *Reichle et al.*, 2002; *Schöniger et al.*, 2012; *Xue and Zhang*, 2014] modeling, etc. As a sequential data assimilation technique, EnKF needs to modify restart files and update model parameters and states simultaneously at each assimilation step, which makes its application inconvenient when the model involves multiple processes [*Emerick and Reynolds*, 2013]. In this situation, computing a global update with all available data is preferred, which leads to the scheme of ensemble smoother (ES) [*Evensen*, 2007; *Van Leeuwen and Evensen*, 1996]. Through only updating model parameters, ES also avoids the inconsistency between updated parameters and states encountered in EnKF. It has been shown that, with much lower computational cost, ES can obtain comparable results as EnKF in some reservoir history matching problems

[*Skjervheim and Evensen*, 2011]. In hydrologic inverse modeling, ES has also found widespread applications, e.g., [*Bailey and Baù*, 2010; *Crestani et al.*, 2013]. However, for strongly nonlinear problems, both EnKF [*Emerick and Reynolds*, 2012; *Gu and Oliver*, 2007; *Lorentzen and Naevdal*, 2011] and ES [*Chen and Oliver*, 2012; *Emerick and Reynolds*, 2013] need some forms of iteration to achieve satisfactory data matches.

As both EnKF and ES rely on the first two statistical moments, they are most suitable for problems with Gaussian distributions. If the distribution of model parameters has multiple modes, using EnKF or ES directly would be problematic. Over the past two decades, there have been several approaches trying to address this issue and extend EnKF or ES to problems with multimodal distributions, most of which are based on cluster analysis. For example, *Elsheikh et al.* [2013] used the K-means algorithm, *Bengtsson et al.* [2003] and some later researchers [*Dovera and Della Rossa*, 2011; *Li et al.*, 2016a; *Smith*, 2007; *Sun et al.*, 2009] used Gaussian mixture models to cluster the samples and update each cluster with EnKF or ES separately. Generally, in these approaches, as we don't know exactly how many modes there are, it would be better to use a relatively large number of clusters. For example, if there are 5 modes (although we don't know this number in advance), setting the number of clusters as 3 would miss some modes and it would be better to set the number of clusters as 5 or a larger number. According to *Elsheikh et al.* [2013], one problem that might be encountered is the stochastic nature of cluster analysis, i.e., different runs of the same inverse algorithm based on cluster analysis may identify different numbers of modes. Moreover, implementing cluster analysis in high-dimensional problems is challenging. Except for adopting cluster analysis, other ways of dealing with multimodal distribution include integrating EnKF with another inverse method, such as particle filter (PF) [*Mandel and Beezley*, 2009], etc.

In this paper, without resorting to the K-means algorithm, Gaussian mixture models, or another inverse algorithm (e.g., PF), we propose a very simple and efficient algorithm, i.e., the iterative local updating ensemble smoother (ILUES), to extend ES to problems with multimodal distributions. For each sample in ES, we define its local ensemble based on an integrated measure of distance to this sample and the measurements. Then we use the scheme of ES to update each local ensemble. In this way, the multimodal distribution of model parameters can be well explored. To achieve satisfactory data matches in strongly nonlinear problems, we adopt an iterative form of ES to assimilate the measurements multiple times.

The remainder of this paper is organized as follows. The detailed formulation of the ILUES algorithm is given in Section 2. To illustrate its performance, five numerical case studies are tested in Section 3. Finally, some conclusions and discussions are provided in Section 4.

## 2. Iterative local updating ensemble smoother

For simplicity, here we represent an arbitrary hydrologic system in the following way:

$$\mathbf{d} = f(\mathbf{m}) + \boldsymbol{\varepsilon}, \tag{1}$$

where $\mathbf{d}$ is a $N_\mathbf{d} \times 1$ vector for the measurements, $f(\cdot)$ is the system model, $\mathbf{m}$ is a $N_\mathbf{m} \times 1$ vector for the uncertain parameters, $\boldsymbol{\varepsilon}$ is a $N_\mathbf{d} \times 1$ vector for the measurement errors. With the noisy measurements $\mathbf{d}$, we can update our knowledge about the unknown model parameters $\mathbf{m}$ via ES:

$$\mathbf{m}_j^a = \mathbf{m}_j^f + \mathbf{C}_{MD}^f (\mathbf{C}_{DD}^f + \mathbf{C}_D)^{-1} [\mathbf{d}_j - f(\mathbf{m}_j^f)], \tag{2}$$

for $j = 1, \ldots, N_e$.

In the above equation, $\mathbf{M}^f = [\mathbf{m}_1^f, \ldots, \mathbf{m}_{N_e}^f]$ is an ensemble of $N_e$ parameter samples randomly drawn from the prior distribution, $\mathbf{M}^a = [\mathbf{m}_1^a, \ldots, \mathbf{m}_{N_e}^a]$ is the updated ensemble conditioned on the measurements $\mathbf{d}$, $\mathbf{C}_{MD}^f$ is the $N_\mathbf{m} \times N_\mathbf{d}$ cross-covariance matrix between $\mathbf{M}^f$ and $\mathbf{D}^f = [f(\mathbf{m}_1^f), \ldots, f(\mathbf{m}_{N_e}^f)]$, $\mathbf{C}_{DD}^f$ is the $N_\mathbf{d} \times N_\mathbf{d}$ auto-covariance matrix of $\mathbf{D}^f$, $\mathbf{C}_D$ is the $N_\mathbf{d} \times N_\mathbf{d}$ covariance matrix of the measurement errors, $\mathbf{d}_j = \mathbf{d} + \boldsymbol{\varepsilon}_j$ is the *j*th realization of the measurements, and $\boldsymbol{\varepsilon}_j$ is a random realization of the measurement errors.

From equation (2), it is obvious that ES only relies on the first two statistical moments. If the prior or the posterior distribution of **m** is multimodal, the direct implementation of ES would be problematic. Nevertheless, being multimodal implies that locally the distribution is still unimodal, which enables the application of ES with a local updating scheme. Based on this idea, we propose a simple and efficient way that identifies and updates $N_e$ local ensembles of $\mathbf{M}^f$ to explore possible multimodal distributions. The local ensemble of the sample $\mathbf{m}_j^f$ ($j = 1, \ldots, N_e$) is identified based on an integrated measure of distance to the measurements **d** and the sample $\mathbf{m}_j^f$:

$$J(\mathbf{m}) = J_1(\mathbf{m})/J_1^{\max} + J_2(\mathbf{m})/J_2^{\max}, \tag{3}$$

where $J_1(\mathbf{m}) = [f(\mathbf{m}) - \mathbf{d}]^T \mathbf{C}_\mathbf{D}^{-1} [f(\mathbf{m}) - \mathbf{d}]$ is the distance between the model responses $f(\mathbf{m})$ and the measurements **d**, and $J_2(\mathbf{m}) = (\mathbf{m} - \mathbf{m}_j^f)^T \mathbf{C}_{\mathbf{MM}}^{-1} (\mathbf{m} - \mathbf{m}_j^f)$ is the distance between the model parameters **m** and the sample $\mathbf{m}_j^f$. Here $\mathbf{C}_{\mathbf{MM}}$ is the $N_\mathbf{m} \times N_\mathbf{m}$ auto-covariance matrix of the model parameters, $J_1^{\max}$ and $J_2^{\max}$ are the maximum values of $J_1(\mathbf{m})$ and $J_2(\mathbf{m})$, respectively. In equation (3), using $J_1^{\max}$ and $J_2^{\max}$ as the scaling factors can make sure that $J_1(\mathbf{m})/J_1^{\max}$ and $J_2(\mathbf{m})/J_2^{\max}$ are within the same range of (0,1], thus neither the $J_1$ part nor the $J_2$ part will dominate.

Then the local ensemble of $\mathbf{m}_j^f$ is the $N_l = \alpha N_e (\alpha \in (0,1])$ samples with the $N_l$ smallest $J$ values, i.e., $\mathbf{M}_j^{l,f} = [\mathbf{m}_{j,1}^f, \ldots, \mathbf{m}_{j,N_l}^f]$. $N_l$ should be large enough so that there are enough samples in the local ensemble to make a reasonable update. Here the factor $\alpha$ represents the ratio between the local ensemble $\mathbf{M}_j^{l,f}$ and the global ensemble $\mathbf{M}^f$. Using the scheme of ES, we can update the corresponding local ensemble:

$$\mathbf{m}_{j,i}^a = \mathbf{m}_{j,i}^f + \mathbf{C}_{\mathbf{MD}}^{l,f}(\mathbf{C}_{\mathbf{DD}}^{l,f} + \mathbf{C}_\mathbf{D})^{-1}[\mathbf{d}_i - f(\mathbf{m}_{j,i}^f)], \tag{4}$$

for $i = 1, \ldots, N_l$. Here $\mathbf{C}_{\mathbf{MD}}^{l,f}$ is the $N_\mathbf{m} \times N_\mathbf{d}$ cross-covariance matrix between $\mathbf{M}_j^{l,f}$ and $\mathbf{D}_j^{l,f} = [f(\mathbf{m}_{j,1}^f), \ldots, f(\mathbf{m}_{j,N_l}^f)]$, $\mathbf{C}_{\mathbf{DD}}^{l,f}$ is the $N_\mathbf{d} \times N_\mathbf{d}$ auto-covariance matrix of $\mathbf{D}_j^{l,f}$, $\mathbf{d}_i = \mathbf{d} + \boldsymbol{\varepsilon}_i$ is the $i$th realization of the measurements. From the updated local ensemble $\mathbf{M}_j^{l,a} = [\mathbf{m}_{j,1}^a, \ldots, \mathbf{m}_{j,N_l}^a]$, we can choose a random sample $\mathbf{m}_j^{l,a}$ as the updated sample of $\mathbf{m}_j^f$ ($j = 1, \ldots, N_e$). In this way, we can well explore the multimodal distribution with the updated global ensemble, $\mathbf{M}^a = [\mathbf{m}_1^{l,a}, \ldots, \mathbf{m}_{N_e}^{l,a}]$.

As stated above, the local ensemble is identified based on an integrated measure of distance both in the space of the model responses ($J_1$) and the model parameters ($J_2$). The role of the $J_1$ part is to filter out the samples that are far away from the posterior region according to the model-data fit, while the role of the $J_2$ part is to filter out the samples that are far away from the mode $\mathcal{M}^*$ that is closest to $\mathbf{m}_j^f$. Through updating the local ensemble of $\mathbf{m}_j^f$, we can obtain the updated parameter sample $\mathbf{m}_j^{l,a}$ that is expected to be close to the mode $\mathcal{M}^*$. With the $N_e$ updated parameter samples in the updated global ensemble $\mathbf{M}^a$, we can identify different modes that may exist in the posterior distribution. If we only use $J_1$ that quantifies the distance between the model responses and the measurement data, we cannot differentiate among different modes and thus cannot solve the multimodal problem. On the other hand, if we only use $J_2$ that quantifies the parametric distance, we can find an ensemble that is close to a certain parameter set. However, it is very likely that the measurement data and the true model parameters are far beyond the bounds of this local ensemble. Then updating this local ensemble is similar to extrapolation and we cannot guarantee to find a good solution. So the $J_1$ part and the $J_2$ part are equally important. By applying equation (3) that considers the two parts simultaneously, we can both differentiate among different modes and make sure that the true state is within or at least not far away from the bounds of the local ensemble. Then updating the local ensembles can provide more robust results.

Another thing that should be noted here is that, when $\alpha = 1$, the local ensemble of a certain sample is actually the entire ensemble, then the updating scheme formulated in equation (4) reduces to that in equation (2) (i.e., the original ES). However, setting $\alpha = 1$ makes the local updating ES unable to handle problems with multiple modes in the posterior. When $\alpha < 1$, the local ensemble of a certain sample is a subset of the global ensemble. Thus the covariance matrices $\mathbf{C}_{\mathbf{MD}}^{f}$ and $\mathbf{C}_{\mathbf{DD}}^{f}$ calculated from the global ensemble will be different from those calculated from the local ensemble (i.e., $\mathbf{C}_{\mathbf{MD}}^{l,f}$ and $\mathbf{C}_{\mathbf{DD}}^{l,f}$). At this point, the local updating ES isn't equivalent to the original ES for the $\alpha < 1$ case. However, the local updating ES is suitable for tackling the multimodal problems, where the performance of the original ES will significantly deteriorate.

In the local updating ES, different local ensembles can share some same samples. At this point, this process is different from the K-means algorithm or Gaussian mixture models. The advantages of this process are twofold. First, as we implement this process with $N_\mathbf{e}$ seeds, it is advantageous in identifying all possible modes when its number is large. Second, if there does not exist any multimodality, different local ensembles would share a considerable number of same samples and thus produce the updated samples that locate in the same mode. In this way, this process can avoid identifying modes erroneously.

For strongly nonlinear problems, an iterative form of ES is usually needed. In this paper, we adopt the simplest one that assimilates the measurements multiple times, which has been integrated into both EnKF [*Emerick and Reynolds*, 2012] and ES [*Emerick and Reynolds*, 2013] for data assimilation in nonlinear problems. At each iteration, we implement the local updating scheme described above on the updated ensemble $\mathbf{M}^a$ obtained from the last iteration. This iterative process is repeated $N_{\text{iter}}$ times. To guarantee that the multiple data assimilation scheme can obtain reasonable results, we need to inflate the covariance matrix of the measurement errors $\mathbf{C_D}$. Here we adopt one simple way that $\mathbf{C_D}$ is multiplied by the predefined iteration number $N_{\text{iter}}$, which has been proven to be able to obtain correct posterior estimates in linear-Gaussian problems using the multiple data assimilation EnKF or ES [*Emerick and Reynolds*, 2012; 2013]. In this scheme, the way to draw realizations of the measurements in equation (4) should be changed accordingly, i.e., $\mathbf{d}_i = \mathbf{d} + \sqrt{N_{\text{iter}}}\mathbf{C_D}^{1/2}\mathbf{r}_{N_\mathbf{d}}$, where $\mathbf{r}_{N_\mathbf{d}} \sim \mathcal{N}(0, I_{N_\mathbf{d}})$. After $N_{\text{iter}}$ iterations, we can obtain good data matches and a converged estimation of the uncertain parameters. Complete scheme of the ILUES algorithm is described in Algorithm 1.

**Algorithm 1** Iterative local updating ensemble smoother

1: Set iteration counter $i = 0$.

2: Generate input ensemble $\mathbf{M}^f = [\mathbf{m}_1^f, ..., \mathbf{m}_{N_e}^f]$ from the prior distribution.

3: Generate output ensemble $\mathbf{D}^f = [f(\mathbf{m}_1^f), ..., f(\mathbf{m}_{N_e}^f)]$ by evaluating the system model.

4: **for** $j = 1, ..., N_e$ **do**

5:      Given $\mathbf{m}_j^f$, calculate the $N_e$ values of $J$ for all samples in $\mathbf{M}^f$ according to equation (3).

6:      Choose the $N_l = \alpha N_e$ samples with the $N_l$ smallest $J$ values as the local ensemble of $\mathbf{m}_j^f$, i.e., $\mathbf{M}_j^{l,f} = [\mathbf{m}_{j,1}^f, ..., \mathbf{m}_{j,N_l}^f]$.

7:      Obtain the updated local ensemble $\mathbf{M}_j^{l,a}$ according to equation (4) with the inflated covariance matrix of the measurement errors and the accordingly generated measurement realizations.

8:      Draw a random sample $\mathbf{m}_j^{l,a}$ from $\mathbf{M}_j^{l,a}$ as the updated sample of $\mathbf{m}_j^f$.

9: **end for**

10: Let $\mathbf{M}^a = [\mathbf{m}_1^{l,a}, ..., \mathbf{m}_{N_e}^{l,a}]$, which is the updated ensemble of $\mathbf{M}^f$.

11: Set $i = i + 1$. If $i = N_{\text{Iter}}$, stop; Otherwise, let $\mathbf{M}^f = \mathbf{M}^a$, go to Step 3.

## 3. Illustrative examples

In this section, we evaluate the ILUES algorithm in five numerical case studies involving nonlinearity and multimodality. The first example is simple and low-dimensional, but it has infinite number of modes in the posterior distribution. This example is used to illustrate the basic ideas of the proposed method. We then test the second example with 100 unknown parameters to demonstrate the performance of ILUES. To show its applicability in complex problems, we further test the ILUES algorithm with three hydrological examples that have multimodal prior distribution, multimodal posterior distribution and a large number ($N_m = 108$) of uncertain parameters, respectively.

3.1. Example 1: A simple case study with infinite number of modes in the posterior

The first example tests the ability of the ILUES algorithm to identify the posterior distribution that has infinite number of modes, which has the following form:

$$y = x_1^2 + x_2^2. \tag{5}$$

In this case, the prior distributions for $x_1$ and $x_2$ are both uniform distributions, $\mathcal{U}(-2, 2)$, the scalar measurement is $d = 1$ with measurement error $\varepsilon \sim \mathcal{N}(0, 0.01^2)$. It is clear that the posterior distribution of the parameters is close to a round circle with radius equal to $\sqrt{d}$, which means that there are infinite number of distinct parameter sets that all can well fit the measurement, i.e., there are infinite number of modes. Although this example is rather simple, it is challenging for the standard ES or the cluster-analysis-based ES to obtain the posterior with infinite number of modes.

Setting the ensemble size $N_e = 400$ and the factor $\alpha = 0.1$ in the ILUES algorithm, the posterior distribution can be well identified within three iterations. The blue dots as shown in Figure 1(a-d) are random samples drawn from the prior distribution and updated samples obtained at the three iterations, respectively. It is clear that the ILUES algorithm is capable of solving inverse problems with infinite number of modes in the posterior distribution. Meanwhile, Figure 1 also demonstrates the necessity of assimilating the measurement multiple times to obtain converged results for nonlinear problems. Here the associated signal to noise ratio defined as the ratio of the average prior root-mean-square error (RMSE) to the average posterior RMSE is 180.17, which indicates a significant reduction of uncertainty in the underlying system.

[Figure 1]

To illustrate the concept of local ensemble, we randomly draw a sample (red diamond) from the prior distribution and plot its local ensemble (black dots) in Figure 1(a). Figure 1(a) indicates that the local ensemble actually locates between the drawn sample and the posterior region, as it is based on an integrated measure of the distance between the model parameters and the drawn sample and the distance between the model response and the measurement. Applying the updating scheme of ES to this local ensemble, we can obtain an updated sample represented by the red diamond in Figure 1(b), which is much closer to the posterior region. The local ensemble of this updated sample is plotted with black dots in Figure 1(b). Similar plots are also shown in Figure 1(c-d).

In the above simulation, the factor $\alpha$ is chosen as 0.1. This factor decides the ratio of the local ensemble over the global ensemble. It is understandable that a smaller $\alpha$ would be more suitable for problems with a large number of modes in the posterior distribution. As we have to make sure that there are enough samples in the local ensemble to make a reasonable update, given a predefined ensemble size $N_e$, $\alpha$ cannot be too small. To illustrate the effect of this factor on the performance of the ILUES algorithm, we test nine different values of $\alpha$ and show the corresponding results in Figure 2 (here $N_e = 400$ with three iterations). In this example, as there are infinite number of modes in the posterior distribution, choosing a large $\alpha$ (e.g., $\alpha > 0.4$) would significantly deteriorate the inversion results. According to our own experience, $\alpha = 0.1$ works well for all our tested examples and thus it is given as the recommended value.

[Figure 2]

Another setting that affects the performance of the ILUES algorithm is the ensemble size $N_e$. As shown in Figure 3 (here $\alpha = 0.1$ with three iterations), when $N_e$ is small (e.g., $N_e = 50$), we will miss a large portion of the posterior region, which greatly underestimates the uncertainty in the model parameters. When $N_e$ is large (e.g., $N_e = 2000$), we can obtain a pretty good result, but it comes with an increased computational cost. Generally speaking, a large $N_e$ is needed for a high-dimensional problem or a problem that has a large number of modes in the posterior. There is a trade-off between the performance and the computational cost when choosing an appropriate $N_e$.

[Figure 3]

Moreover, the updated sample of $\mathbf{m}_j^f$, i.e., $\mathbf{m}_j^{l,a}$, is randomly drawn from the updated local ensemble $\mathbf{M}_j^{l,a}$. If not choosing randomly, but selecting the updated sample that has the smallest $J$ value seems to be an appealing option. As shown in Figure 4, given different settings of the ensemble size $N_e$ and the factor $\alpha$, selecting the "best" sample (solid lines) would always give better data matches than choosing a random sample (dashed lines). Here the $y$ axis in Figure 4 is for the log-transformed RMSE (Log RMSE) between the simulated model outputs and the measurement averaged over the ensemble. However, as we will demonstrate in the following example, this option may cause biased inversion results in some specific problems. In Figure 4, it is again shown that smaller values of $\alpha$ can usually give better data matches. As we should preserve enough samples in each local ensemble to make a reasonable update via ES, a small $\alpha$ should come with a relatively large ensemble size $N_e$, i.e., a relatively high computational cost.

[Figure 4]

In equation (3), the measure for the local ensemble of $\mathbf{m}_j^f$ ($j = 1, \ldots, N_e$) assigns equal weights to the normalized distance between the model responses and the measurements $\mathbf{d}$ (i.e., $J_1(\mathbf{m})/J_1^{\max}$) and the normalized distance between the model parameters $\mathbf{m}$ and the sample $\mathbf{m}_j^f$ (i.e., $J_2(\mathbf{m})/J_2^{\max}$). Here we can also assign different weights to the two normalized distances:

$$J(\mathbf{m}) = J_1(\mathbf{m})/J_1^{\max} + b \cdot J_2(\mathbf{m})/J_2^{\max}, \qquad (6)$$

where $b \in (0, \infty)$. In Figure 5, we systematically study the effect of the factors $\alpha$ and $b$ on the performance of the ILUES algorithm (here $N_e = 400$ with three iterations and the "random" option). When $\alpha < 0.1$, smaller values of $b$ (e.g., 0.1 and 0.01) can obtain better data matches. However, when $\alpha \geq 0.1$, it is better to choose a relatively large value of $b$ (e.g., $b > 0.1$). This is because choosing a big value of $b$ would make the local ensemble relatively close to the sample $\mathbf{m}_j^f$ and relatively far away from the measurements $\mathbf{d}$. If $\alpha$ is very small (e.g., $\alpha = 0.01$), the local ensemble would have a very small size and it may miss the samples that are close to the measurements $\mathbf{d}$, which could cause dissatisfactory data matches. When $\alpha$ is relatively large, the local ensemble can keep some samples that are close to the measurements. In this situation, preserving the local properties of the sample $\mathbf{m}_j^f$ might matter more. In our experience, $b = 1$ could provide more robust results than other values. Moreover, in many papers working on inverse problems, similar objective functions to equation (3) have been formulated (although the scaling factors $J_1^{\max}$ and $J_2^{\max}$ might not be used) [*Chen and Oliver*, 2012; *Zhou et al.*, 2014], where the contribution of the parametric distance has the same weight as the distance in the model responses, i.e., $b = 1$. Thus, $b = 1$ is used as the default value in the following examples.

[Figure 5]

3.2. Example 2: A 100-dimensional case study with multimodal posterior

To show the performance of the ILUES algorithm in problems with more unknown model parameters, we test the second example:

$$y = x_1^2 + x_2^2 + \cdots + x_{100}^2. \tag{7}$$

Here the prior distributions for $x_1 \sim x_{99}$ are $\mathcal{U}(0,1)$, and for $x_{100}$ is $\mathcal{U}(-10,10)$. The scalar measurement in this case is $d = 87.68$, with measurement error $\varepsilon \sim \mathcal{N}(0, 1^2)$. It is expected that the posterior distribution of $x_{100}$ is bimodal, i.e., using either $x_{100}$ or $-x_{100}$ we can obtain the same model response when other model parameters are the same.

[Figure 6]

As this problem has 100 unknown model parameters, a relatively large ensemble size is chosen in the ILUES algorithm. In this case, $N_\mathrm{e} = 1000$ and $\alpha = 0.1$ are used. From Figure 6 we can find that, within five iterations, the bimodality of $x_{100}$ in the posterior distribution can be well identified by the ILUES algorithm. As shown in Figure 7(a), although the simulated model outputs from the prior samples have a large uncertainty level, they can converge to the measurement within five iterations. In this case study, the ratio of the average prior RMSE to the average posterior RMSE is 20.84. Besides, the RMSE between $f(\mathbf{m}_{\text{true}})$ and the actual measurement $d$ is 1.16, the mean of the posterior RMSE is 0.91, and the 95% confidence interval of the posterior RMSE is [0.029 2.69], which is close to the results of MCMC simulation (mean: 0.94, 95% confidence interval:[0.038 2.88]). Here we have to admit that using only 1000 samples is far from enough to fully characterize the 100-dimensional posterior distribution. However, it is still a good way to make an accurate prediction of the system.

[Figure 7]

In this example, we also test the option that selects the updated sample that has the smallest $J$ value from the updated local ensemble in the ILUES algorithm. However, as shown in Figure 7(b), it will cause a biased inversion result that has an abnormally large variance. From Figure 7(b) we can also find that in the last three iterations, many samples would stay near where they were at the last iteration. It may be because the "best" sample in each local ensemble is usually closest to the "true" state, which would receive the smallest update. Moreover, this update would become even smaller at later iterations, which could prevent a sufficient update of the model parameters. Thus, this option is not very robust. In the following examples, we will choose the updated sample randomly from the updated local ensemble and this setting will not be further specified.

3.3. Example 3: A rainfall-runoff model with multimodal prior

The third example tests the ability of the ILUES algorithm to deal with problems whose prior distributions have multiple modes. Here we consider a more practical case, which is based on a widely used rainfall-runoff model, HYMOD, developed by *Boyle* [2000]. This model connects a simple rainfall-excess model [*Moore*, 1985] to a series of linear slow and quick reservoirs within a watershed. There are five uncertain parameters in HYMOD, i.e., the maximum water storage capacity of the watershed,

$C_{max}[L]$, the degree of spatial variability of soil moisture capacity, $b_{exp}[-]$, the distribution factor for the flow between the slow and the quick reservoirs, $\beta[-]$, the residence time of the slow reservoirs, $R_s[T]$ and the residence time of the quick reservoirs, $R_q[T]$, respectively. This example is included in the DREAM software package developed by *Vrugt* [2016] and it is modified and used in this case study. Here the prior distributions for $C_{max}$ and $b_{exp}$ are multimodal and represented by Gaussian mixture models, i.e., $p(C_{max}) = 1/3\mathcal{N}(100, 20^2) + 1/3\mathcal{N}(250, 20^2) + 1/3\mathcal{N}(400, 20^2)$ and $p(b_{exp}) = 1/3\mathcal{N}(0.5, 0.1^2) + 1/3\mathcal{N}(1, 0.1^2) + 1/3\mathcal{N}(1.5, 0.1^2)$, respectively. While the prior distributions for $\beta$, $R_s$ and $R_q$ are uniform distributions, whose ranges are listed in Table 1. The stream flow measurements are generated from one set of true model parameters $\mathbf{m}_{true}$ as listed in Table 1 with additive measurement errors $\boldsymbol{\varepsilon} \sim \mathcal{N}(\mathbf{0}, \boldsymbol{\sigma}^2)$, where $\boldsymbol{\sigma} = 0.1 \times f(\mathbf{m}_{true})$.

[Table 1]

Choosing the ensemble size $N_e = 300$ and the factor $\alpha = 0.1$, the ILUES algorithm can accurately estimate the model parameters within five iterations, as shown in Figure 8. Here the ratio of the average prior RMSE to the average posterior RMSE is 5.30. Besides, the RMSE between $f(\mathbf{m}_{true})$ and the actual measurements $\mathbf{d}$ is 5.35, the mean of the posterior RMSE is 5.37, and the 95% confidence interval of the posterior RMSE is [5.32 5.46], which is close to the results of MCMC simulation (mean:5.35, 95% confidence interval: [5.30 5.42] ). Compared with example 1, although there are more parameters in this example, a smaller $N_e$ is capable of quantifying the parametric uncertainties, as the number of modes is not large. However, there is still a trade-off between the performance and the computational cost. If a very small $N_e$ is chosen, there will be a risk of obtaining biased inversion results.

[Figure 8]

3.4. Example 4: Contaminant source identification with multimodal posterior

In this example, we consider a contaminant source identification problem in steady-state saturated groundwater flow. As shown in Figure 9, the $20[L] \times 10[L]$ domain has constant-head conditions at the left ($12[L]$) and right ($11[L]$) boundaries, no-flow conditions at the lower and upper boundaries, respectively. The conductivity and porosity of the aquifer are homogeneous, whose values are known as $K = 8[LT^{-1}]$ and $\theta = 0.25[-]$, respectively. Then we can obtain a uniform background flow from left to right. In this flow field, some amount of contaminant is released from a point source. The contaminant source is characterized by five parameters, i.e., $\mathbf{m} = [x_s, y_s, S_s, t_{on}, t_{off}]$, which means that the contaminant is released at $(x_s, y_s)[L]$ from time $t_{on}[T]$ to $t_{off}[T]$ with a constant mass-loading rate $S_s[MT^{-1}]$. The prior distributions for the five parameters are uniform, whose ranges are listed in Table 2. To infer these parameters, concentration measurements are collected from a single well denoted by the blue circle in Figure 9 at $t = [6, 8, 10, 12, 14][T]$ with measurement errors $\varepsilon \sim \mathcal{N}(0, 0.01^2)$. The true values of the model parameters $\mathbf{m}_{true}$ that generate the measurements are also listed in Table 2.

[Figure 9]

[Table 2]

The governing equations for the steady-state saturated groundwater flow are:

$$\frac{\partial}{\partial x_i}\left(K_i \frac{\partial h}{\partial x_i}\right) = 0, \tag{8}$$

and

$$v_i = -\frac{K_i}{\theta}\frac{\partial h}{\partial x_i}, \tag{9}$$

where $h[L]$ represents hydraulic head, $K_i[LT^{-1}]$ and $v_i[LT^{-1}]$ represent hydraulic conductivity and pore water velocity along the respective coordinate axis $x_i[L](i = 1,2)$, respectively.

The advection dispersion equation for the contaminant transport is:

$$\frac{\partial(\theta C)}{\partial t} = \frac{\partial}{\partial x_i}\left(\theta D_{ij}\frac{\partial C}{\partial x_j}\right) - \frac{\partial}{\partial x_i}(\theta v_i C) + q_s C_s, \tag{10}$$

where $C[ML^{-3}]$ represents molar concentration of the dissolved contaminant; $t[T]$ is time; $q_s[T^{-1}]$ and $C_s[ML^{-3}]$ represent flow rate per unit volume of aquifer and concentration of the contaminant source; $D_{ij}[L^2T^{-1}]$ represents hydrodynamic dispersion tensor, whose principal components ($D_{xx}$ and $D_{yy}$) and cross terms ($D_{xy}$ and $D_{yx}$) are defined as:

$$\begin{cases} D_{xx} = (\alpha_L v_x^2 + \alpha_T v_y^2)/|\boldsymbol{v}|, \\ D_{yy} = (\alpha_L v_y^2 + \alpha_T v_x^2)/|\boldsymbol{v}|, \\ D_{xy} = D_{yx} = (\alpha_L - \alpha_T)v_x v_y/|\boldsymbol{v}|, \end{cases} \quad (11)$$

where $\alpha_L$ and $\alpha_T$ represent longitudinal and transverse dispersivities, $v_x$ and $v_y$ represent components of the pore water velocity $\boldsymbol{v}$ along $x$ and $y$ directions, $|\boldsymbol{v}|$ is the magnitude of $\boldsymbol{v}$, respectively. Here the longitudinal and transverse dispersivities are known as $\alpha_L = 0.3[L]$ and $\alpha_T = 0.03[L]$, respectively. The governing equations for the groundwater flow and solute transport are numerically solved with MODFLOW [*Harbaugh et al.*, 2000] and MT3DMS [*Zheng and Wang*, 1999], respectively.

[Figure 10]

To estimate the model parameters, we implement the ILUES algorithm with $N_e = 400$ and $\alpha = 0.1$. As shown in Figure 10, the posterior distribution of $y_s$ is bimodal. In this case study, the ratio of the average prior RMSE to the average posterior RMSE is 67.21. Besides, the RMSE between $f(\mathbf{m}_{\text{true}})$ and the actual measurements $\mathbf{d}$ is 0.0069, the mean of the posterior RMSE is 0.0086, and the 95% confidence interval of the posterior RMSE is $[0.0048\ 0.014]$, which is close to the results of MCMC simulation (mean:0.0081, 95% confidence interval:[0.0045 0.014]). To verify that the inversion result obtained by the ILUES algorithm is reasonable, we also show the parameter estimation results obtained by the MCMC simulation. In this case, the DREAM algorithm developed by Vrugt is adopted, whose efficiency has been shown in inverse problems with multimodal distributions [*Vrugt*, 2016]. Here we use eight parallel chains in the DREAM algorithm, each of which has a length of 2000, i.e., the total number of model evaluations is 16,000. The Gaussian likelihood function is used to evaluate the goodness-of-fit between the model outputs and the measurements. As shown in Figure 11, the trace plots of the model parameters obtained by DREAM are very similar to those obtained by the ILUES algorithm.

[Figure 11]

Here we also implement ES that assimilates the measurements multiple times to estimate the model parameters. As shown in Figure 12, using the same ensemble size and the same number of iterations, ES with multiple data assimilation cannot accurately characterize the bimodal posterior distribution of $y_s$, although it can still reduce the uncertainties of $x_s$, $t_{on}$ and $t_{off}$ whose posterior distributions are unimodal.

[Figure 12]

3.5. Example 5: Contaminant source identification with 108 unknown parameters

To demonstrate the performance of the ILUES algorithm in inverse problems with many more unknown model parameters, we further test a more complex contaminant source identification problem. In this example, instead of considering a source with a constant strength, we consider a time-varying source strength, which is characterized by 6 parameters in 6 time segments, i.e., $S_{si}[MT^{-1}]$ during $i:i+1[T]$, for $i=1,\dots,6$. Therefore, along with the source location $(x_s, y_s)$, there are 8 parameters that characterize the contaminant source. Again, these parameters are assumed to follow uniform distributions, whose ranges are listed in Table 3.

[Table 3]

In this example, we consider the heterogeneity of the conductivity field whose log-transformed values $Y = \log(K)$ at two arbitrary locations $(x_1, y_1)$ and $(x_2, y_2)$ are assumed to be correlated in the following form:

$$C_Y(x_1, y_1; x_2, y_2) = \sigma_Y^2 \exp\left(-\frac{|x_1 - x_2|}{\lambda_x} - \frac{|y_1 - y_2|}{\lambda_y}\right), \qquad (12)$$

where $\sigma_Y^2 = 1$ is the variance, $\lambda_x = 10[L]$ and $\lambda_y = 5[L]$ are the correlation lengths along *x* and *y* directions, respectively. Here we use the Karhunen-Loève (KL) expansion [*Zhang and Lu*, 2004] to parameterize the log-transformed conductivity field:

$$Y(\mathbf{x}) \approx \bar{Y}(\mathbf{x}) + \sum_{i=1}^{N_{\text{KL}}} \sqrt{\tau_i} s_i(\mathbf{x}) \xi_i, \tag{13}$$

where $\bar{Y}(\mathbf{x}) = 2$ is the mean component, $\tau_i$ and $s_i(\mathbf{x})$ are eigenvalues and eigenfunctions of the correlation function described in equation (12), $\xi_i (i = 1, \dots, N_{\text{KL}})$ are independent standard Gaussian random variables. In this case, 100 KL terms are kept, i.e., $N_{\text{KL}} = 100$, which can preserve about 94.7% of the field variance, i.e., $\sum_{i=1}^{100} \tau_i / \sum_{i=1}^{\infty} \tau_i \approx 94.7\%$.

Thus, there are 108 unknown model parameters in this case, i.e., the 8 parameters for the contaminant source and the 100 KL terms for the log-transformed conductivity field. To infer these parameters, we collect concentration measurements at $t = [4, 5, 6, 7, 8, 9, 10, 11, 12][T]$ and hydraulic head measurements at the fifteen wells denoted by the blue squares in Figure 9. The measurement errors for the concentration and hydraulic head are all assumed to be independent and Gaussian with zero means and standard deviations of $0.005[ML^{-3}]$ and $0.005[L]$, respectively. The reference log-transformed conductivity field and true values of the contaminant source parameters are shown in Figure 14(a) and Table 3, respectively.

[Figure 13]

In this case with 108 unknown model parameters, a large ensemble size $N_e = 2000$ is chosen in the ILUES algorithm with $\alpha = 0.1$. As shown in Figure 13, the contaminant source parameters can be accurately identified within seven iterations. Meanwhile, three realizations, the mean and variance of the posterior log conductivity field are presented in Figure 14(b-f), which clearly demonstrate the estimation accuracy of the log-transformed conductivity field. Here the ratio of the average prior RMSE to the average posterior RMSE of the concentration data is 586.47, and for the head data the ratio is 27.44. Besides, the RMSE between $f(\mathbf{m}_{\text{true}})$ and the measurements $\mathbf{d}$, the means and 95% confidence intervals of the posterior RMSEs obtained by ILUES and MCMC are listed in Table 4. It is shown that, the RMSE value between $f(\mathbf{m}_{\text{true}})$ and the measurements $\mathbf{d}$ slightly deviates from the confidence intervals of the posterior RMSEs obtained by both ILUES and MCMC. Thus, the relatively large number of unknown model parameters ($N_\mathbf{m} = 108$) pose a challenge for both algorithms in accurate uncertainty quantification.

[Figure 14]

[Table 4]

In this case, there isn't any parameter whose distribution is obviously multimodal, but we can still use similar settings as those used in inverse problems with multimodal distributions. At this point, the ILUES algorithm has an advantage over previous cluster analysis-based methods, which need to make a subtle choice of the number of clusters in advance. On the other hand, the ILUES algorithm usually needs much fewer model evaluations than MCMC. In this example with 108 unknown model parameters, even the state-of-the-art DREAM algorithm would need hundreds of thousands of model evaluations.

It should be noted here that for a high-dimensional problem (e.g., $N_\mathbf{m} > 100$), the ensemble size of a few thousand might not be enough to fully quantify the parametric uncertainty. When the input-output relationship of the high-dimensional problem is complex and nonlinear, we had better set a larger ensemble size and more iterations. Moreover, in complicated high-dimensional problems, two samples that have close values of $J$ as defined in equation (3) may not be actually similar. In this situation, the ability of the ILUES algorithm in identifying multiple posterior modes may compromise.

## 4. Conclusions and discussions

In this paper, to extend the ensemble smoother (ES) to inverse problems with multimodal distributions, we propose a simple and efficient algorithm, i.e., the iterative local updating ensemble smoother (ILUES). For each sample in ES, we define the local ensemble based on an integrated measure of the distance between the model responses and the measurements and the distance between the model parameters and the originally drawn sample. Then we use the scheme of ES to update each local ensemble. In this way, the multimodal distribution can be well explored. To achieve satisfactory data matches in nonlinear problems, a simple iterative form of ES that assimilates the measurements multiple times is adopted.

Five numerical case studies are tested to show the performance of the proposed method. The first example demonstrates the ability of the ILUES algorithm to tackle posterior distribution with infinite number of modes. In this simple case study, we systematically illustrate the basic ideas of the proposed method. The second example is similar to the first one, but has many more unknown parameters ($N_\mathbf{m} = 100$). The other three case studies are inverse problems in hydrologic modeling, which consider possible multiple modes in the prior or posterior distributions. All these case studies successfully show the performance of the proposed method in adequately quantifying parametric uncertainties of complex systems, no matter whether the multimodal distribution exists.

In the above examples, we only consider the measurement error. While in many situations, the model structural error should also be considered. In that case, one has to explicitly express $\boldsymbol{\varepsilon}$ as the measurement error plus the model structural error, i.e., $\boldsymbol{\varepsilon}_{\text{total}} = \boldsymbol{\varepsilon}_{\text{measurement}} + \boldsymbol{\varepsilon}_{\text{model}}$. As the distribution of the structural error is usually unknown, we can estimate the parameters that describe the error distribution together with the unknown model parameters in the ILUES algorithm. Similar strategies have been applied in parameter estimation problems with MCMC, e.g., [*Vrugt et al.*, 2009a]. In another approach, the model structural error can be quantified with a data-driven approach (e.g., Gaussian process [*Xu and Valocchi*, 2015] ) during the model calibration period. When multiple model proposals are available, we can adopt the framework of Bayesian model averaging to rigorously consider the model structural uncertainty [*Rojas et al.*, 2008; *Ye et al.*, 2004], which has also been applied in the framework of ensemble Kalman filter, e.g., [*Xue and Zhang*, 2014].

In this paper, the multimodality stems from the system nonlinearity and scarcity of measurement data. In reservoir simulation, the multimodality originated from strongly non-Gaussian parameter field (e.g., multi-facies and channelized permeability fields) is also drawing people's attention [*Jafarpour and Mclaughlin*, 2009]. In this situation, people have to adopt additional techniques, e.g., the level set method [*Chang*

*et al.*, 2010] and normal-score transform [*Zhou et al.*, 2011] to handle the discretely distributed parameter fields. These issues will be addressed in our future works.

# Acknowledgments

Computer codes and data used are available at:

https://www.researchgate.net/publication/315805609_MATLAB_code_of_the_ILUES_algorithm

The authors would like to thank the three anonymous reviewers and the editors for their constructive comments on an earlier version, which led to a much improved paper.

This work is supported by the National Natural Science Foundation of China (grants 41721001 and 41771254). Guang Lin would like to acknowledge the support of NSF grant DMS-1555072 and the U.S. Department of Energy, Office of Science, Office of Advanced Scientific Computing Research, Applied Mathematics program as part of the Multifaceted Mathematics for Complex Energy Systems ($M^2ACS$) project and part of the Collaboratory on Mathematics for Mesoscopic Modeling of Materials project.

We would also like to acknowledge Jasper Vrugt from UC Irvine for providing the codes of DREAM, an efficient MCMC algorithm.

# Tables

Table 1 Prior ranges and true values of model parameters in the third example

| Parameter | $C_{\max}[L]$ | $b_{\exp}[-]$ | $\beta[-]$ | $R_s[T]$ | $R_q[T]$ |
|---|---|---|---|---|---|
| Range | [1 500] | [0.1 2] | [0.1 0.99] | [0 0.1] | [0.1 0.99] |
| True value | 417.416 | 1.464 | 0.362 | 0.0254 | 0.694 |

Table 2 Prior ranges and true values of model parameters in the fourth example

| Parameter | $x_s[L]$ | $y_s[L]$ | $S_s[MT^{-1}]$ | $t_{\text{on}}[T]$ | $t_{\text{off}}[T]$ |
|---|---|---|---|---|---|
| Range | [3 5] | [3 7] | [10 13] | [3 5] | [9 11] |
| True value | 3.854 | 5.999 | 11.044 | 4.897 | 9.075 |

Table 3 Prior ranges and true values of contaminant source parameters in the fifth example

| Parameter | Range | True value |
|---|---|---|
| $x_s[L]$ | [3 5] | 3.520 |
| $y_s[L]$ | [4 6] | 4.437 |
| $S_{s1}[MT^{-1}]$ | [0 8] | 5.692 |
| $S_{s2}[MT^{-1}]$ | [0 8] | 7.883 |
| $S_{s3}[MT^{-1}]$ | [0 8] | 6.306 |
| $S_{s4}[MT^{-1}]$ | [0 8] | 1.485 |
| $S_{s5}[MT^{-1}]$ | [0 8] | 6.872 |
| $S_{s6}[MT^{-1}]$ | [0 8] | 5.552 |

Table 4 The RMSE between $f(\mathbf{m}_{\text{true}})$ and the measurements **d**, the means and 95% confidence intervals of the posterior RMSEs obtained by ILUES and MCMC

| | Measurements | ILUES | | MCMC | |
|---|---|---|---|---|---|
| | | mean | 95% interval | mean | 95% interval |
| Concentration $[ML^{-3}]$ | 0.0045 | 0.0069 | [0.0057 0.0089] | 0.0114 | [0.0105 0.0121] |
| Head $[L]$ | 0.0033 | 0.0044 | [0.0034 0.0058] | 0.0104 | [0.0078 0.0135] |

# Figures

Figure 1. (a) Random samples drawn from the prior distribution and (b-d) updated samples obtained at the three iterations. The local ensemble (black dots) of the sample denoted by the red diamond is shown in each subplot.

Figure 2. With different values of the factor $\alpha$, the obtained results of parameter estimation. Here $N_e = 400$ with three iterations.

Figure 3. With different values of the ensemble size $N_e$, the obtained results of parameter estimation. Here $\alpha = 0.1$ with three iterations.

Figure 4. With different values of the ensemble size $N_e$ and the factor $\alpha$, the log-transformed RMSE between the simulated model outputs and the measurement averaged over the ensemble. Here the dashed lines are for the scenario that we randomly choose the updated sample from the updated local ensemble and the solid lines are for the scenario that we select the updated sample with the smallest $J$ value from the updated local ensemble.

Figure 5. With different values of $\alpha$ and $b$, the log-transformed RMSE between the simulated model outputs and the measurement averaged over the ensemble. Here $N_e = 400$ with three iterations and the "random" option.

Figure 6. Trace plot of $x_{100}$ obtained by the ILUES algorithm in the second example. Here $N_e = 1000$ and $\alpha = 0.1$.

Figure 7. Simulated model outputs of the ILUES algorithm by (a) choosing the updated sample randomly from the updated local ensemble, and (b) selecting the updated sample with the smallest $J$ value from the updated local ensemble.

Figure 8. Trace plots of model parameters obtained by the ILUES algorithm in the third example. Here $N_e = 300$ and $\alpha = 0.1$.

Figure 9. Flow domain for the fourth and fifth examples. The potential area of the contaminant source is represented by the red dashed rectangle. The measurement locations for the fourth and fifth examples are denoted by the blue circle and the blue squares, respectively.

Figure 10. Trace plots of model parameters obtained by the ILUES algorithm in the fourth example. Here $N_e = 400$ and $\alpha = 0.1$.

Figure 11. Trace plots of model parameters obtained by the DREAM algorithm in the fourth example.

Figure 12. Trace plots of model parameters obtained by ES with multiple data assimilation in the fourth example.

Figure 13. Trace plots of contaminant source parameters obtained by the ILUES algorithm in the fifth example. Here $N_e = 2000$ and $\alpha = 0.1$.

Figure 14. (a) Reference log-transformed conductivity field, (b-d) three posterior realizations of the log-transformed conductivity field, (e) mean estimate of the log-transformed conductivity field, and (f) estimation variance of the log-transformed conductivity field.

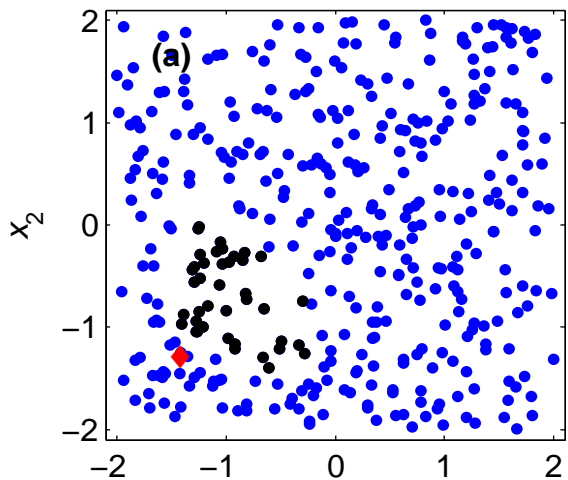
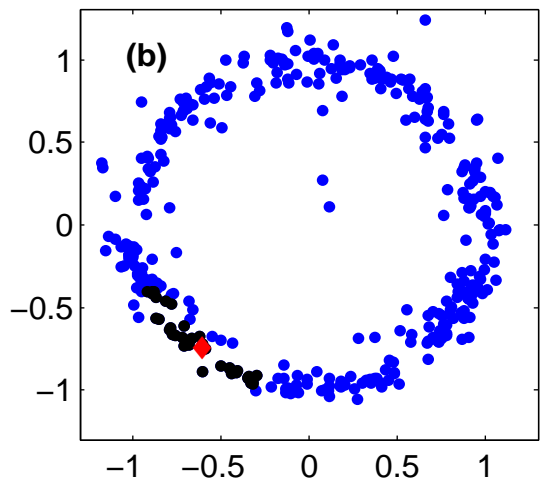
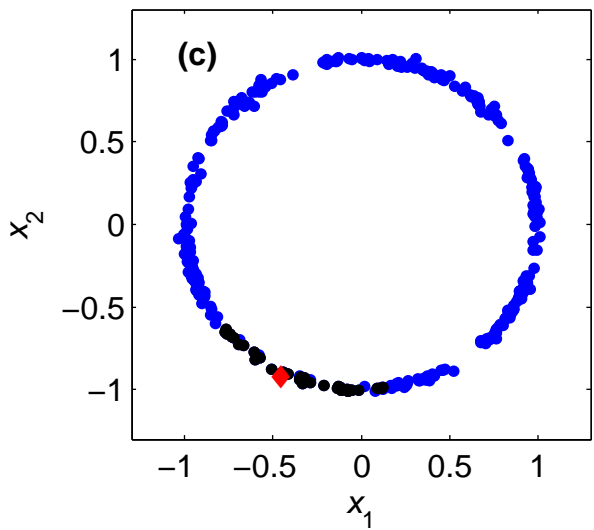
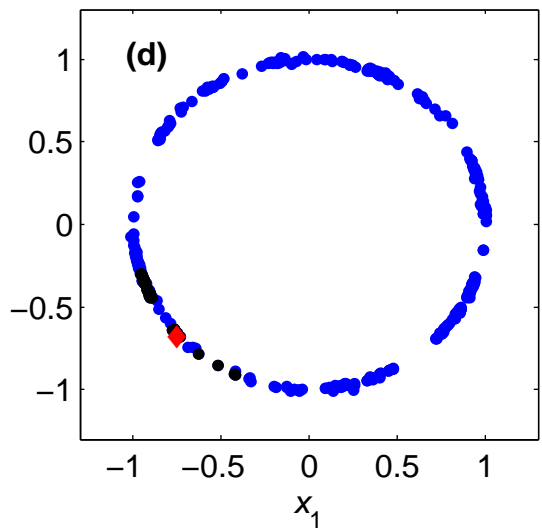

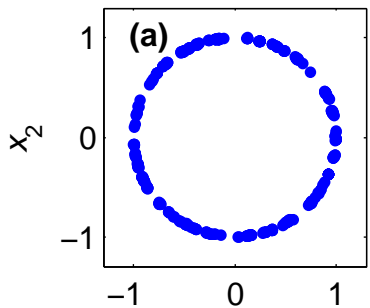 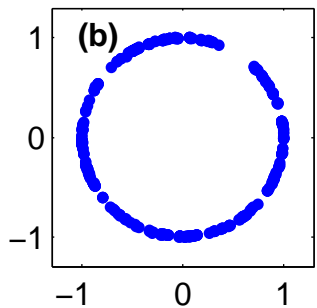 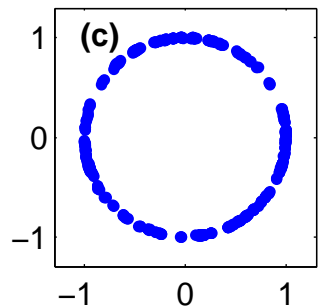
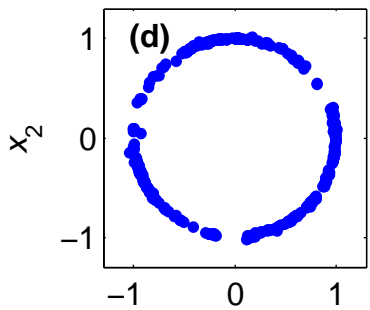 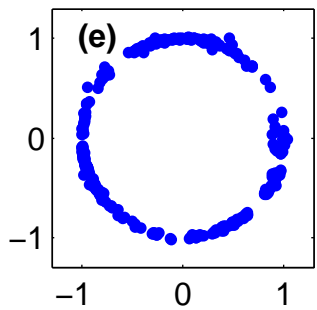 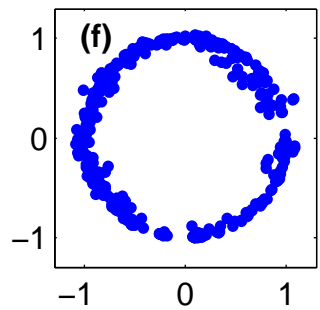
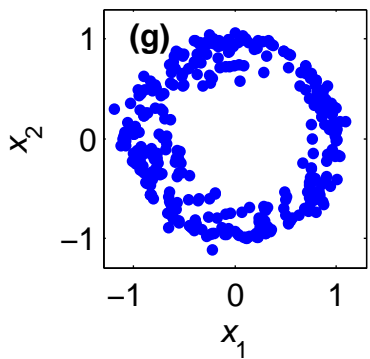 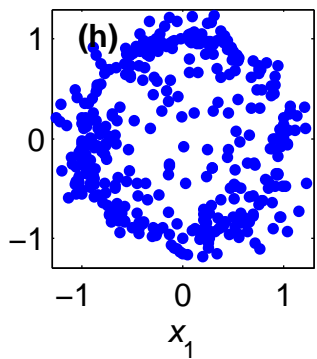 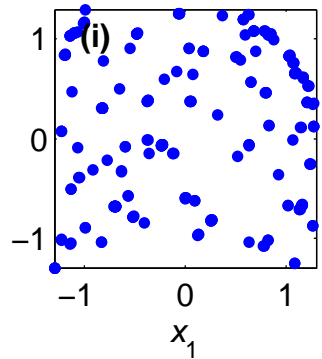

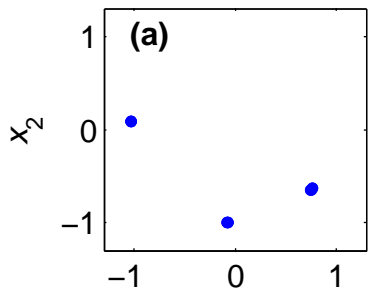 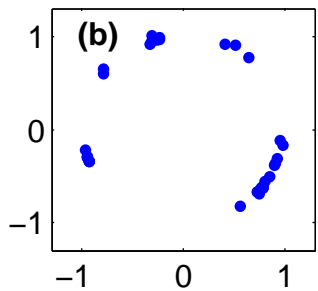 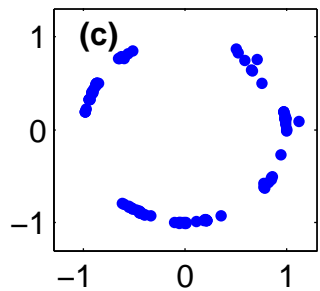
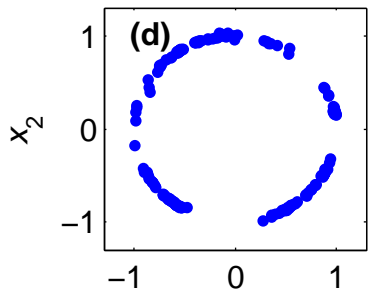 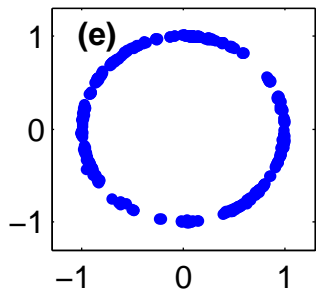 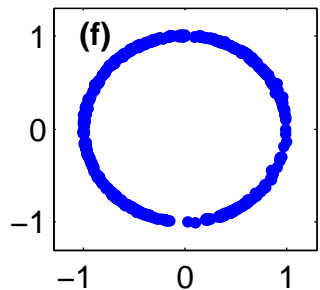
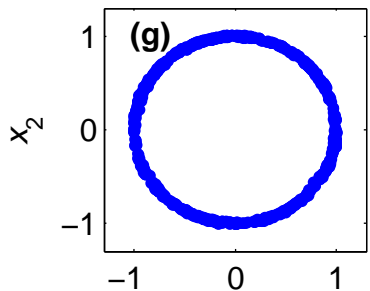 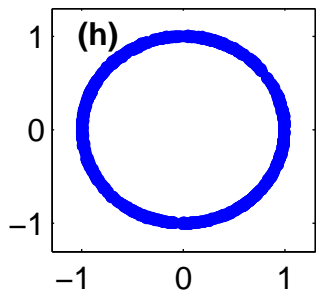 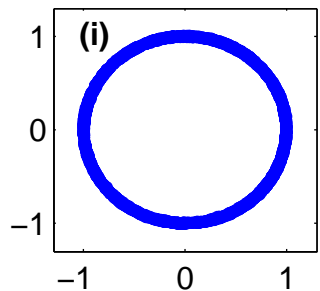

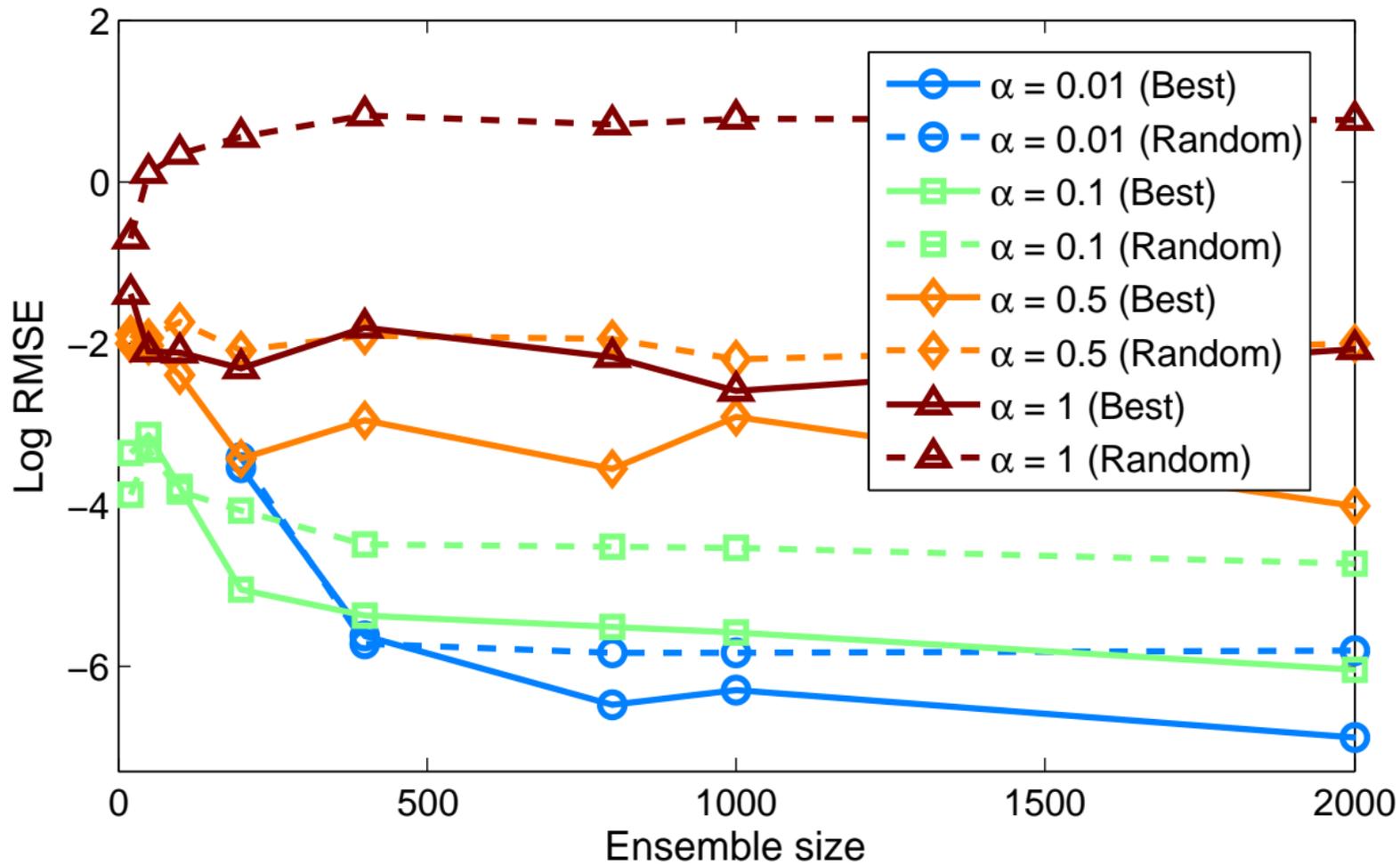

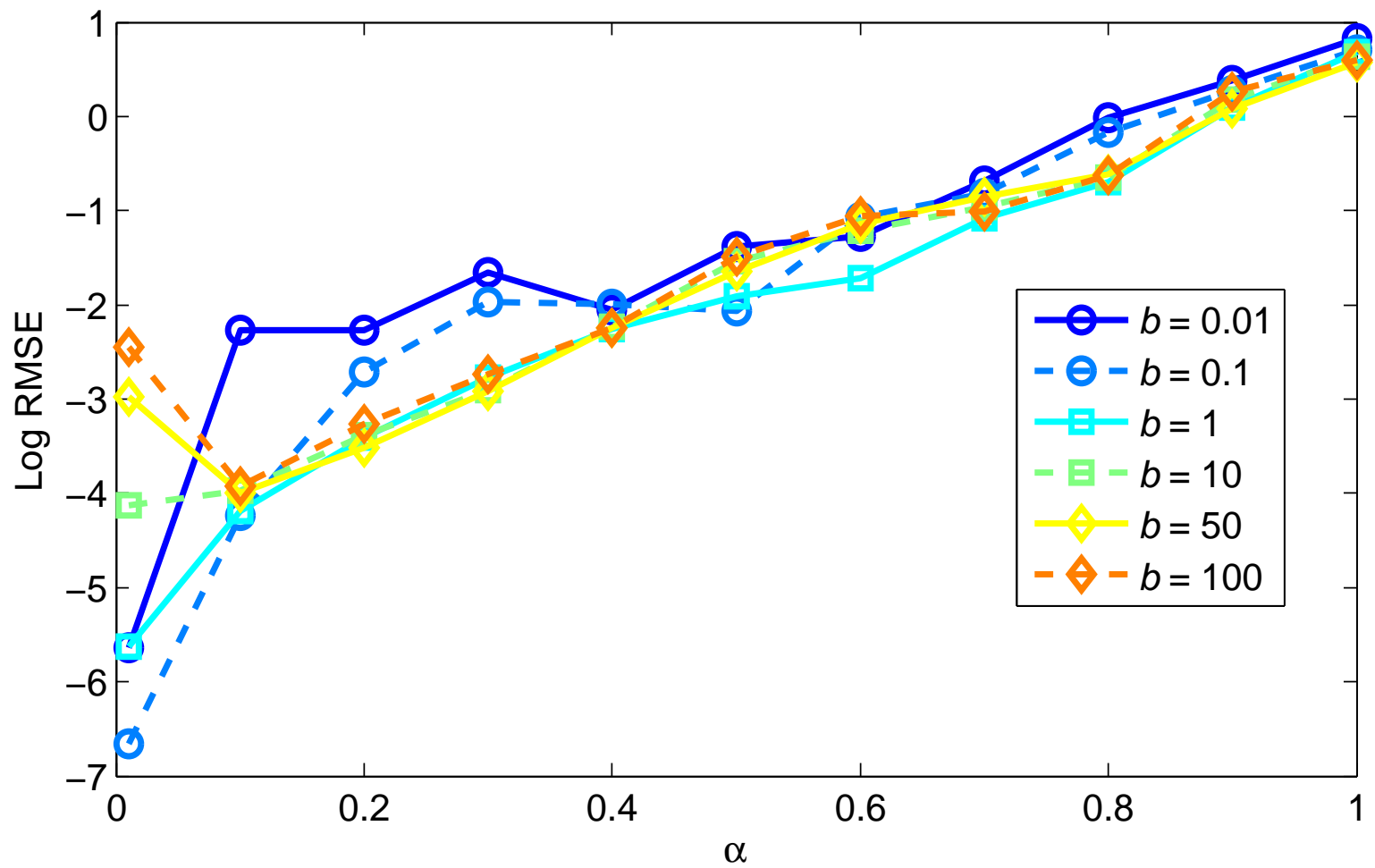

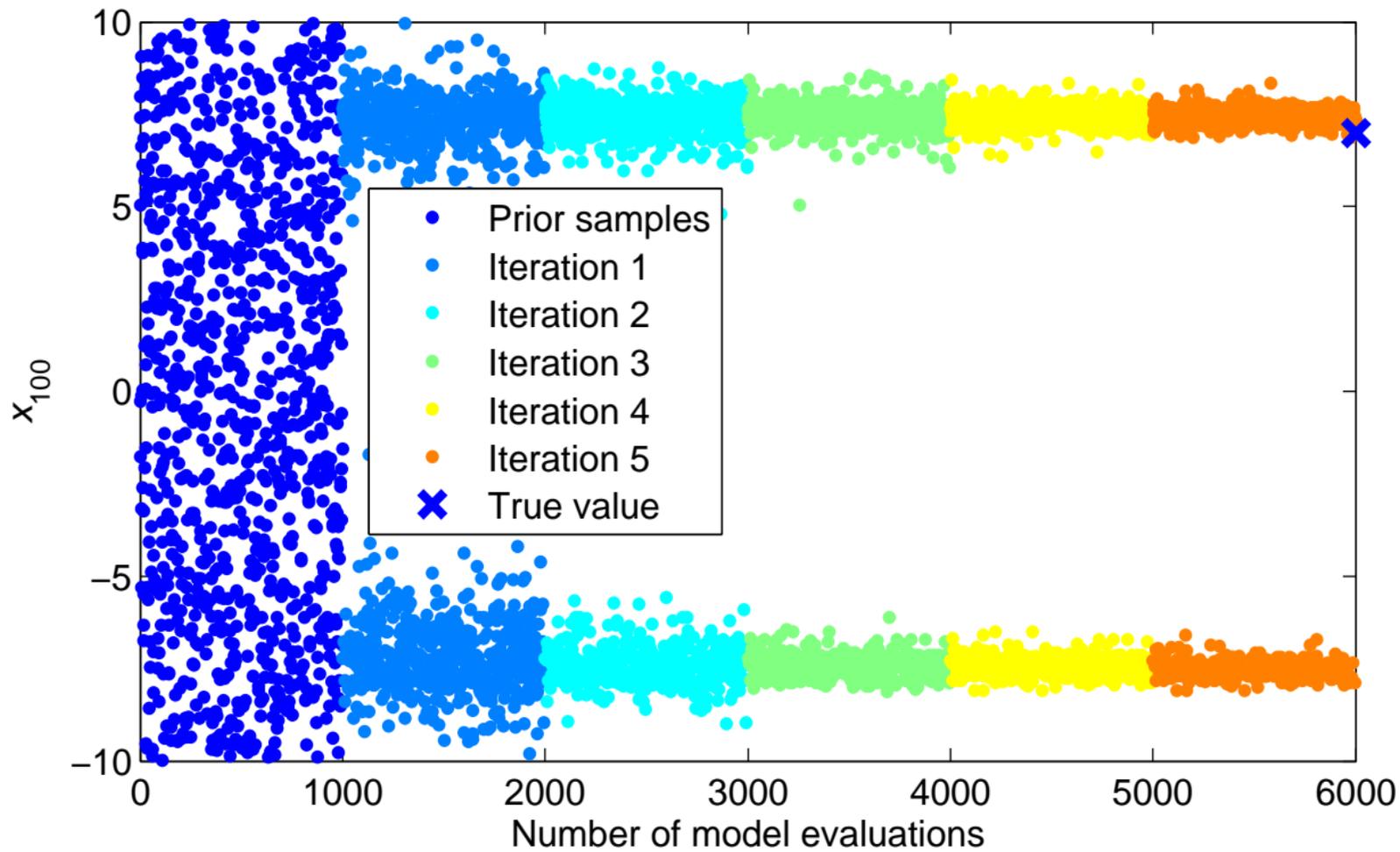

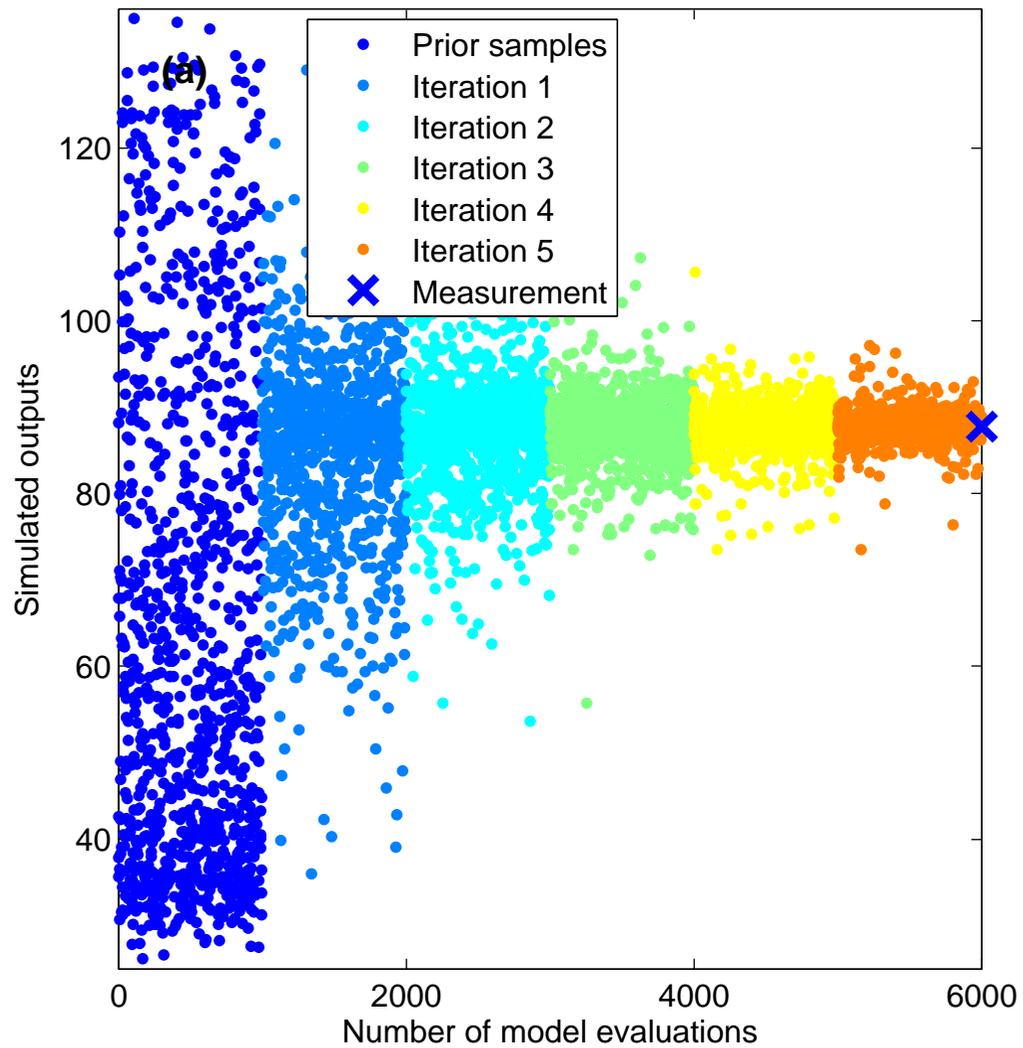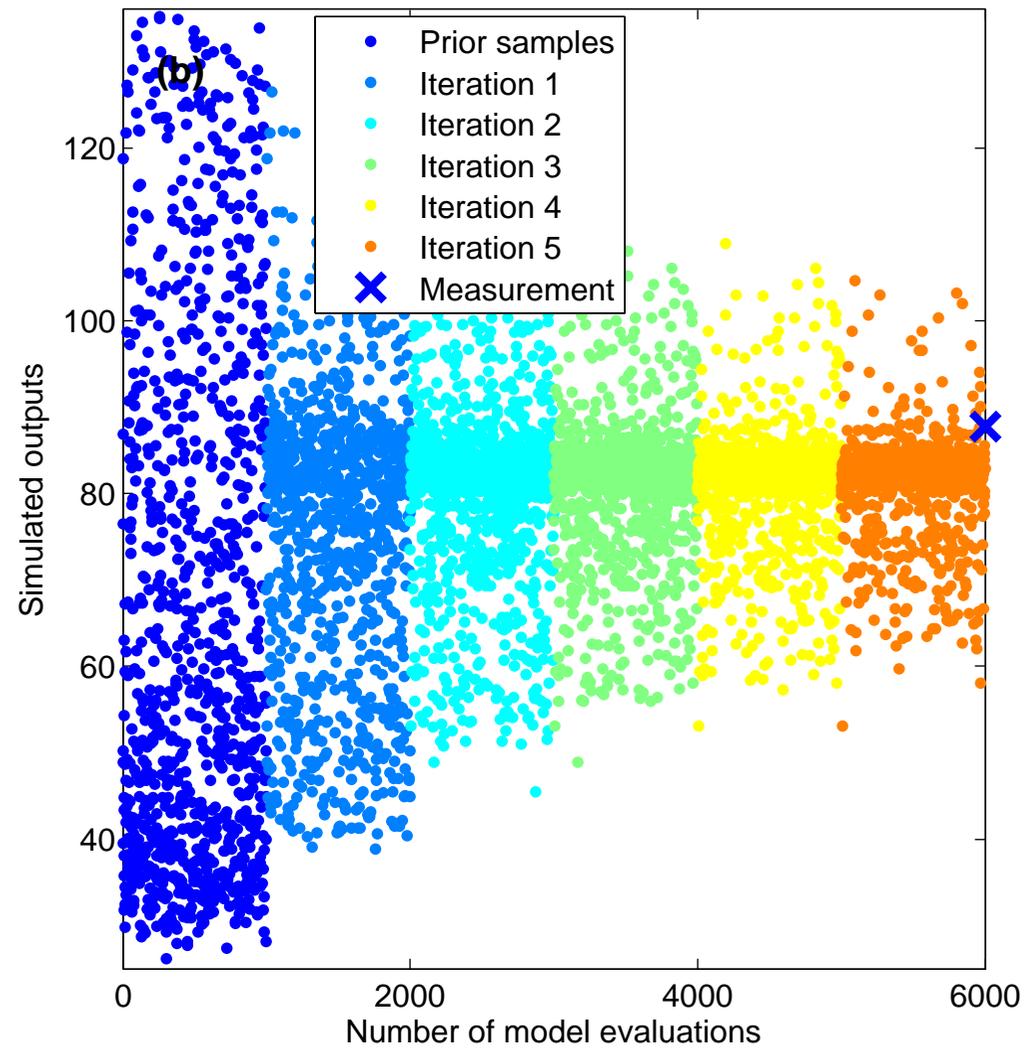

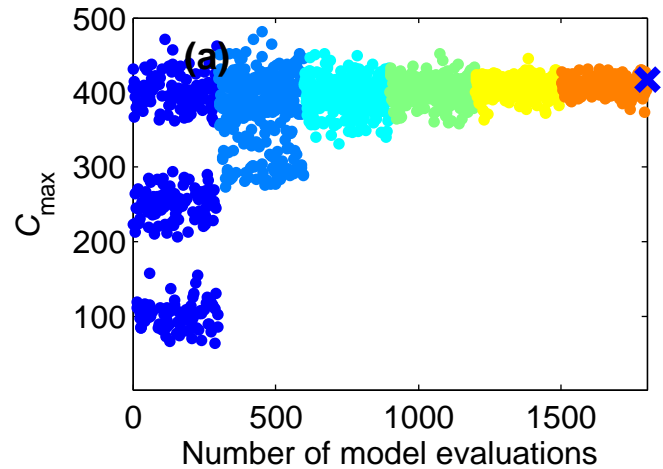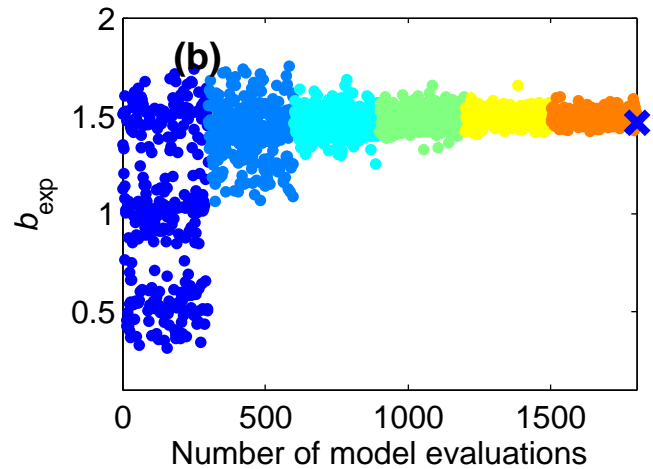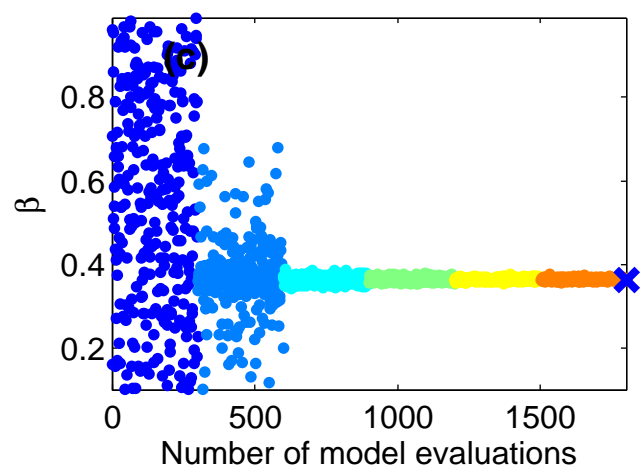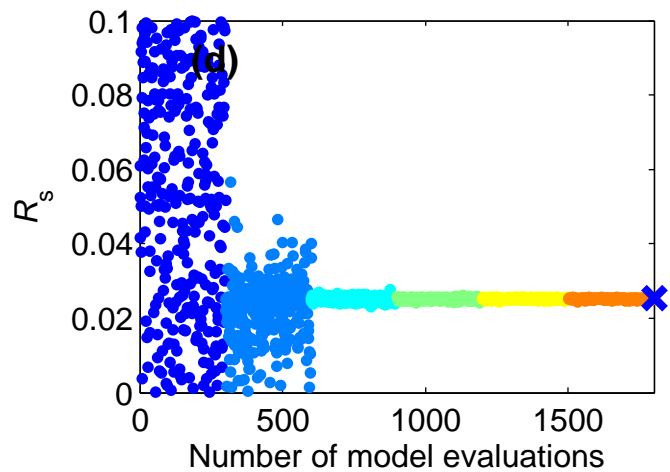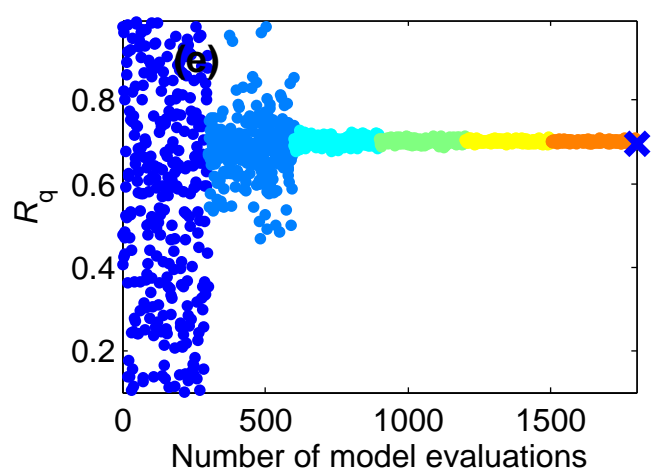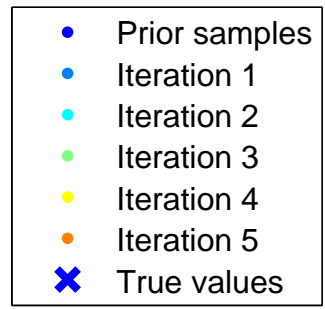

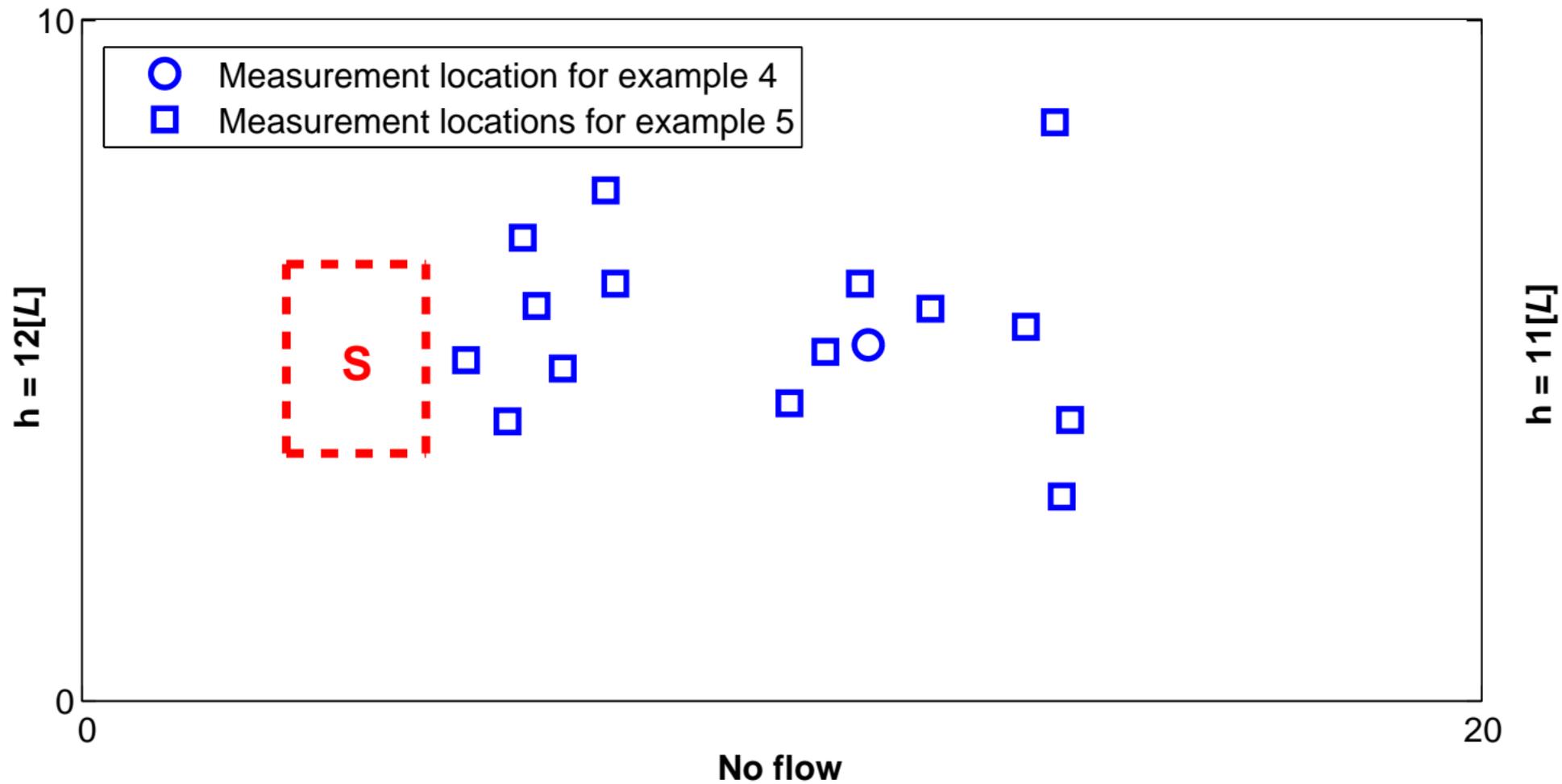

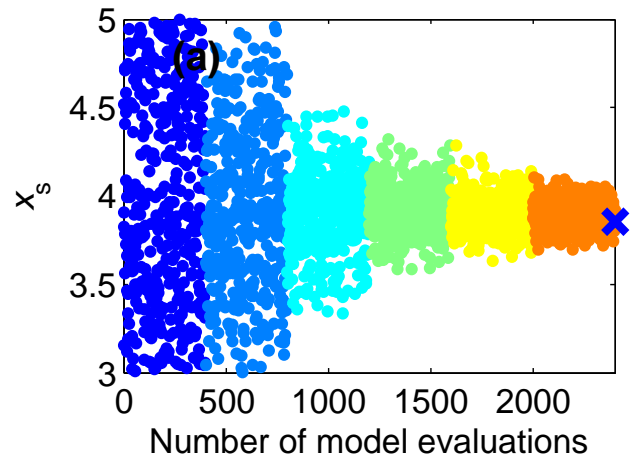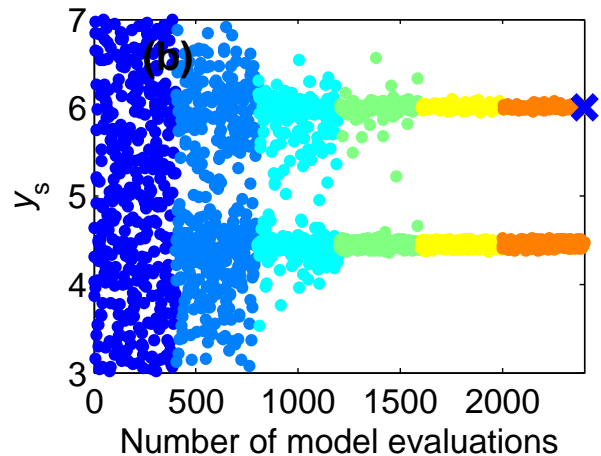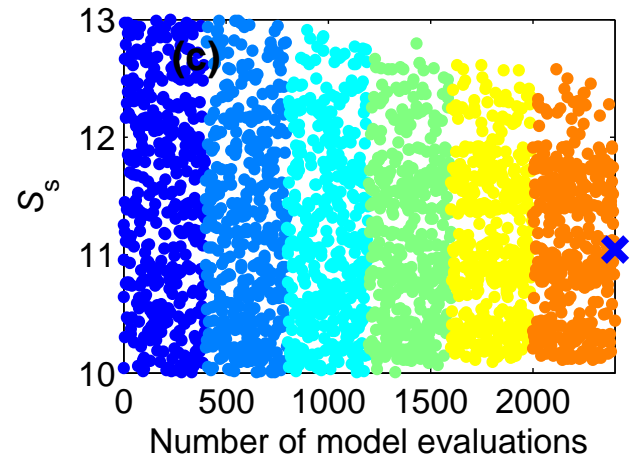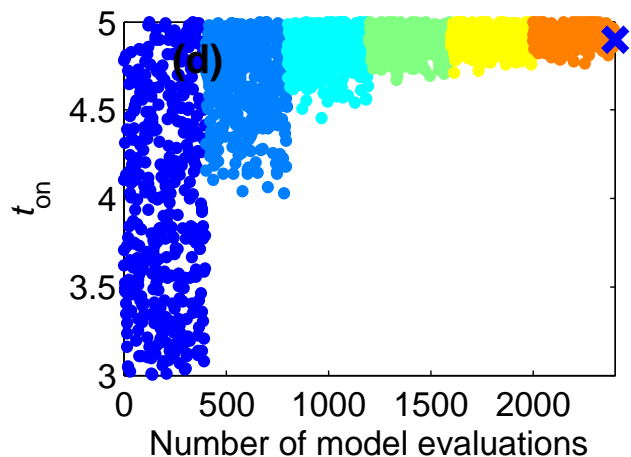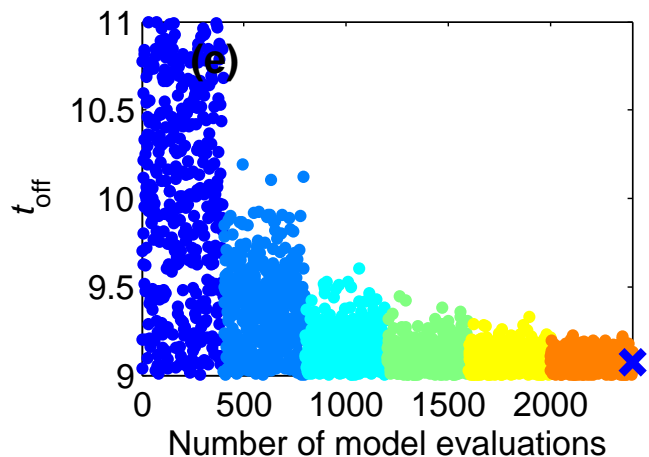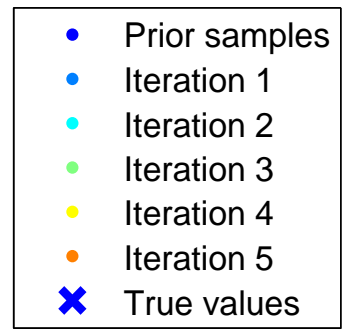

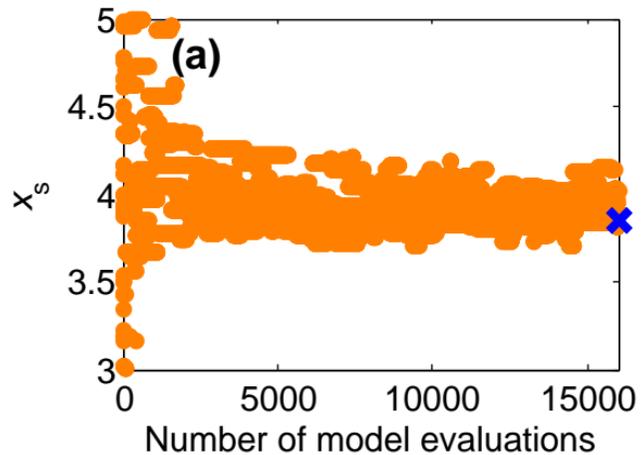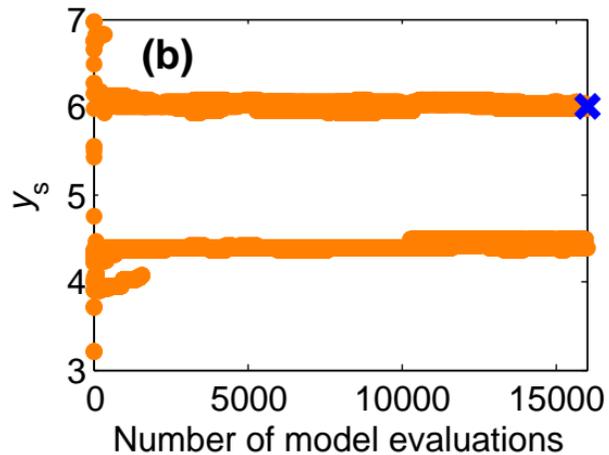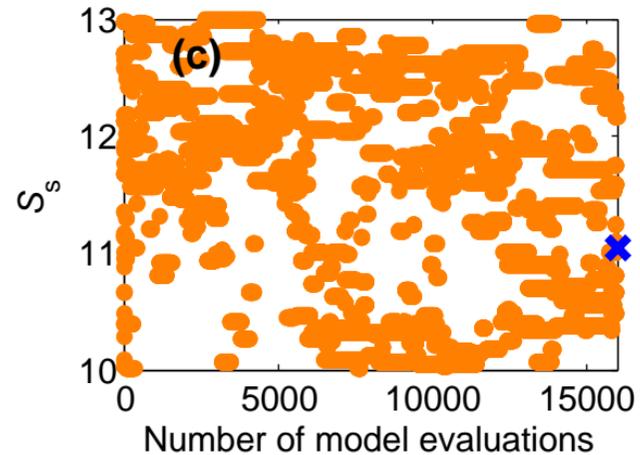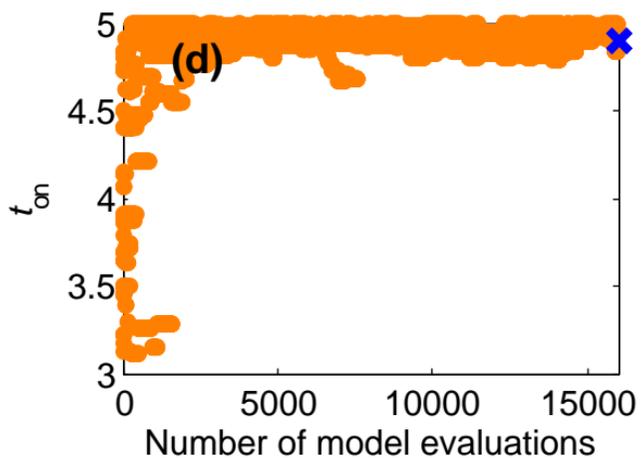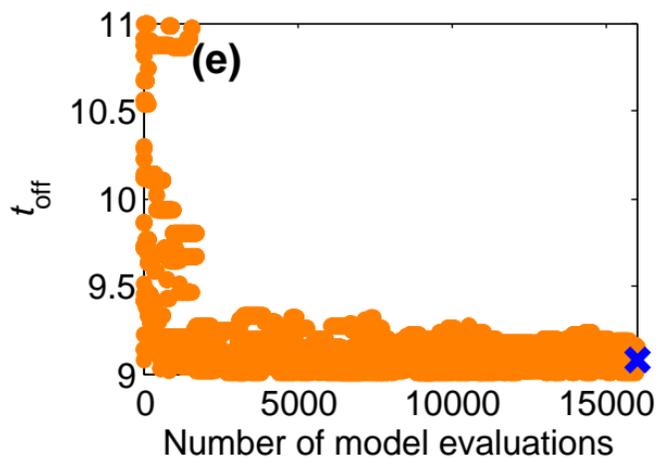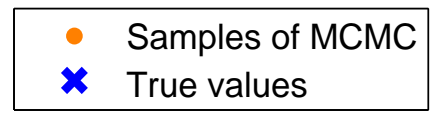

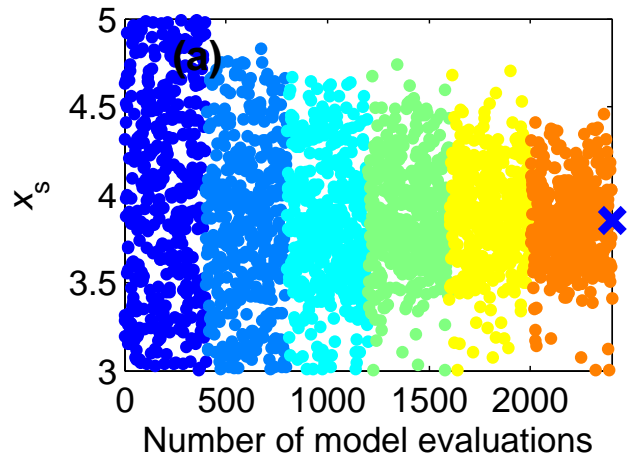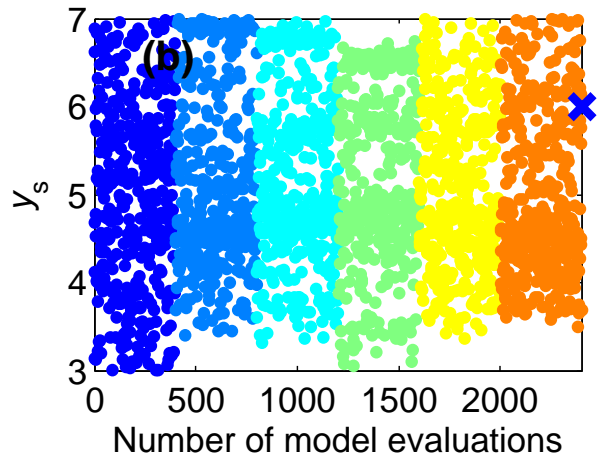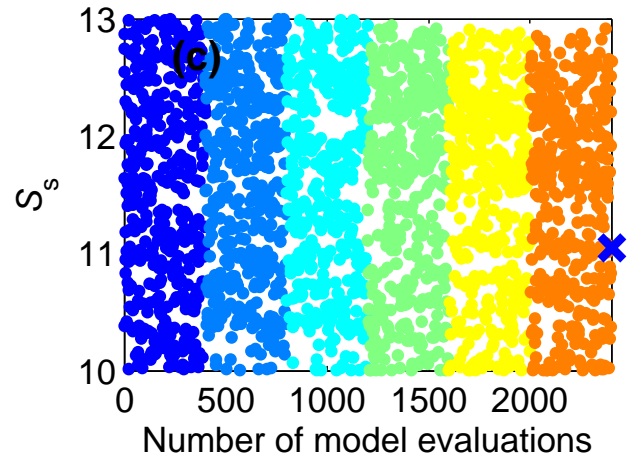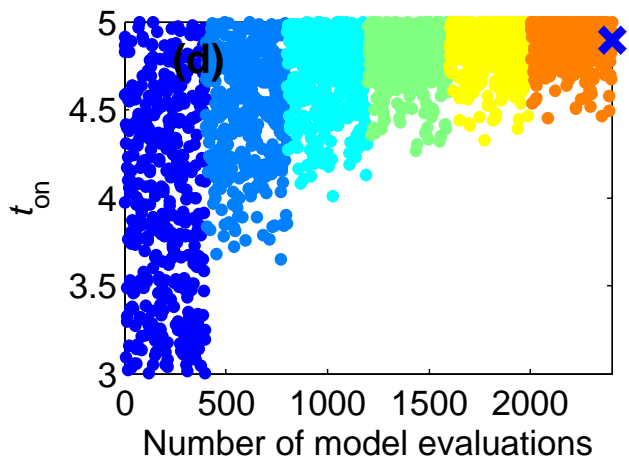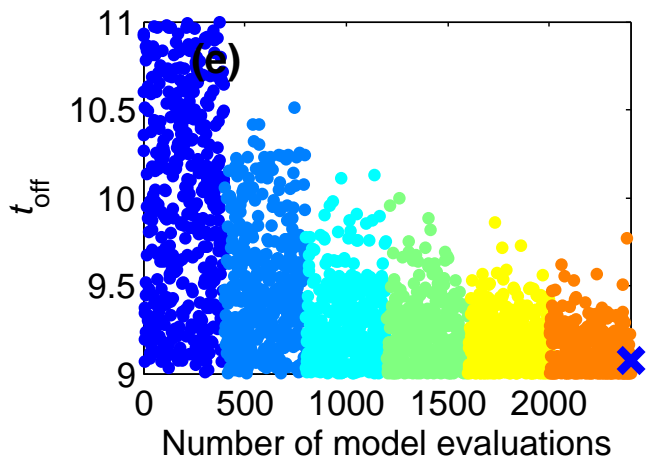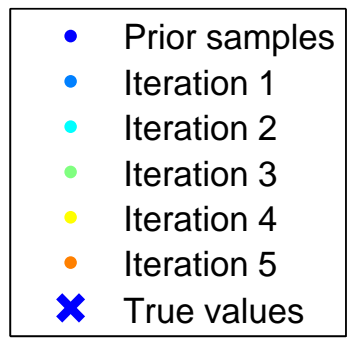

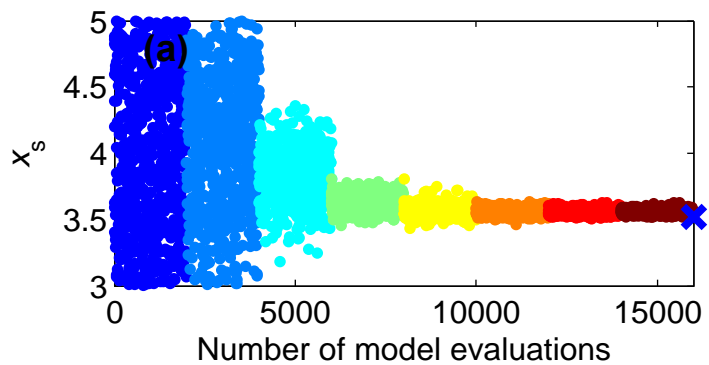
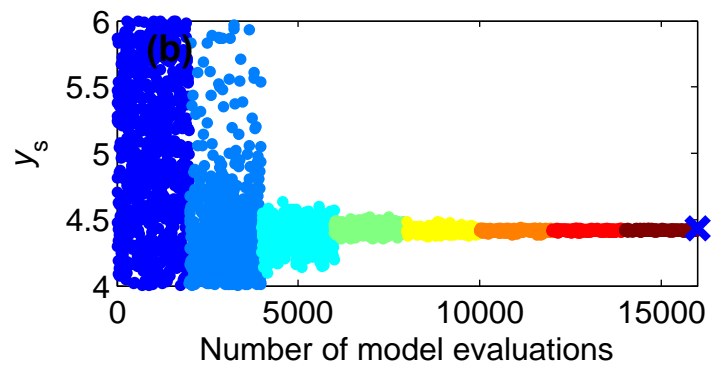
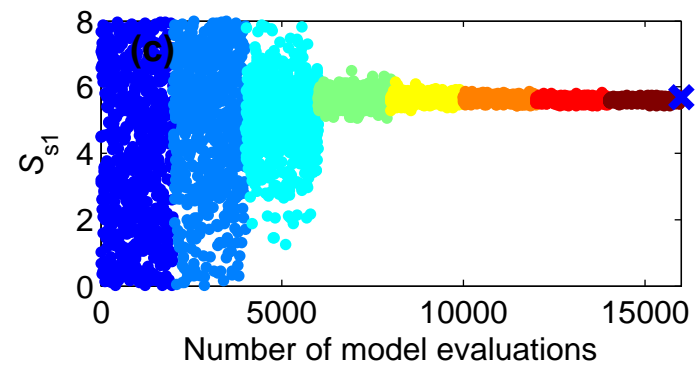
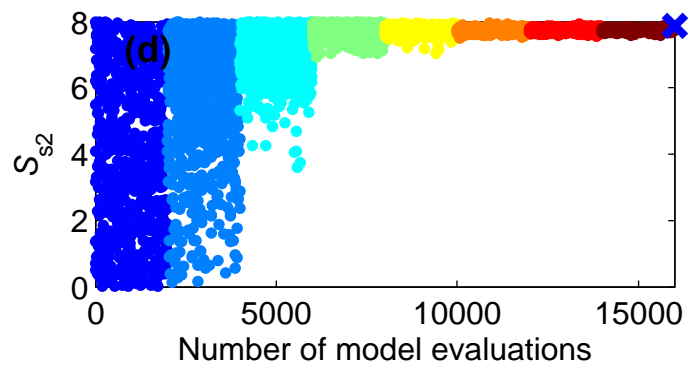
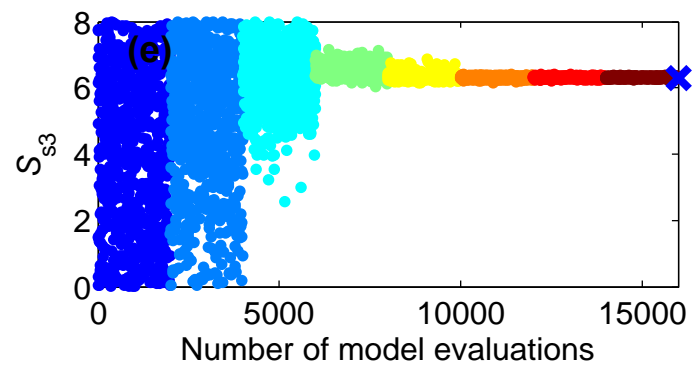
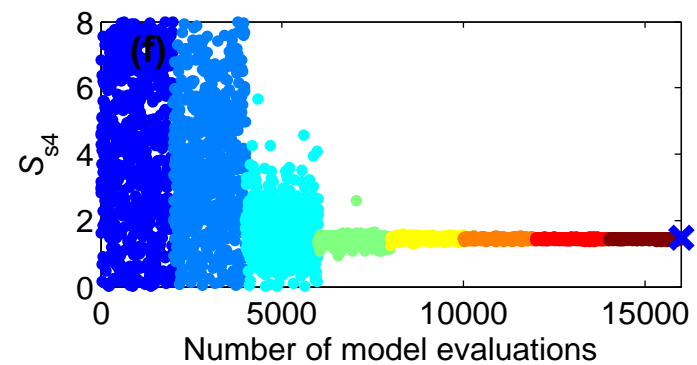
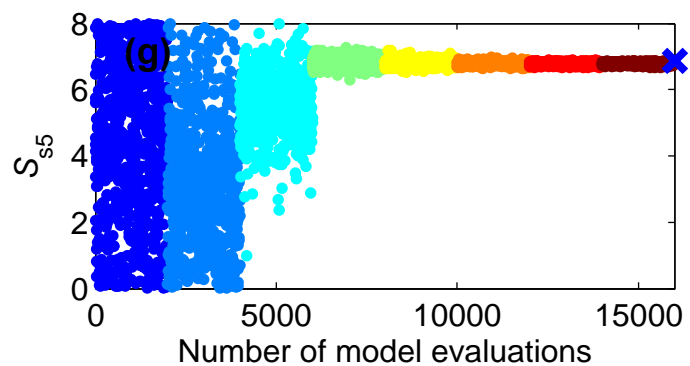
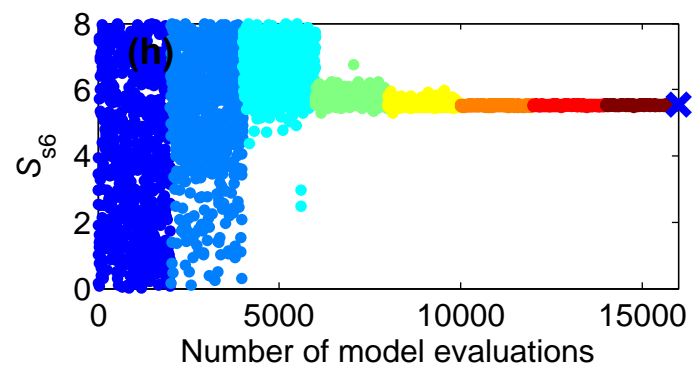
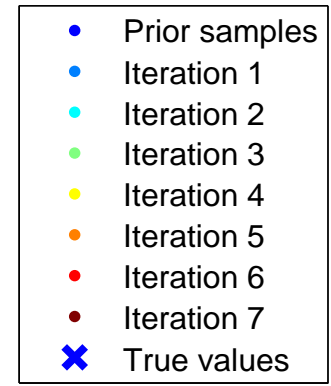

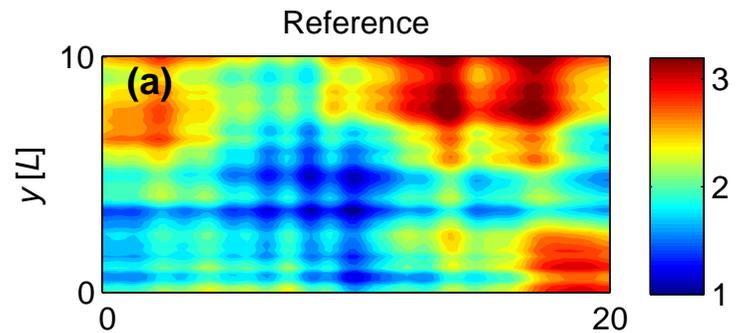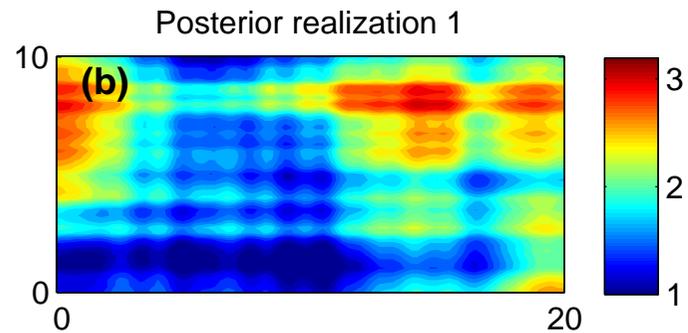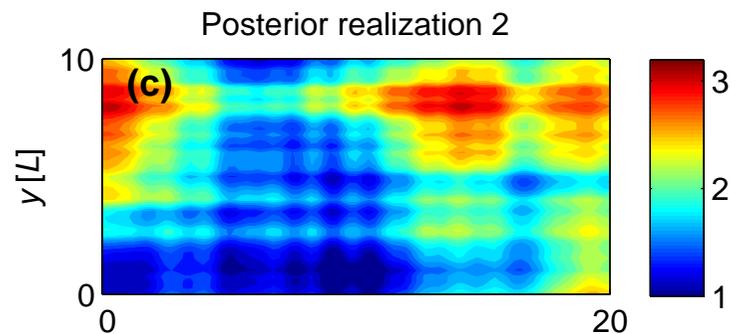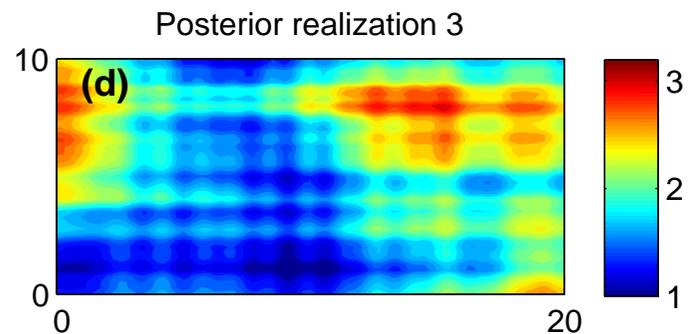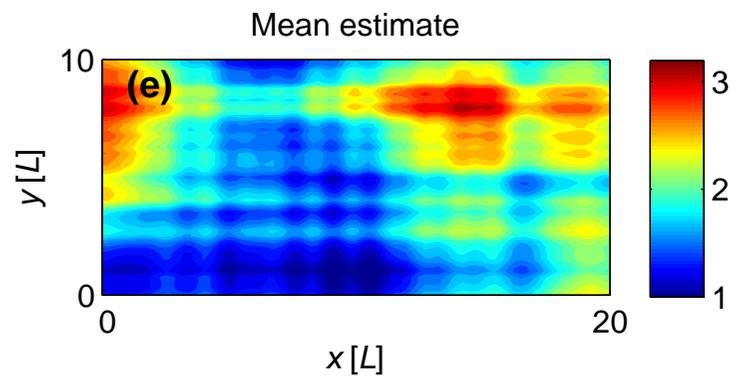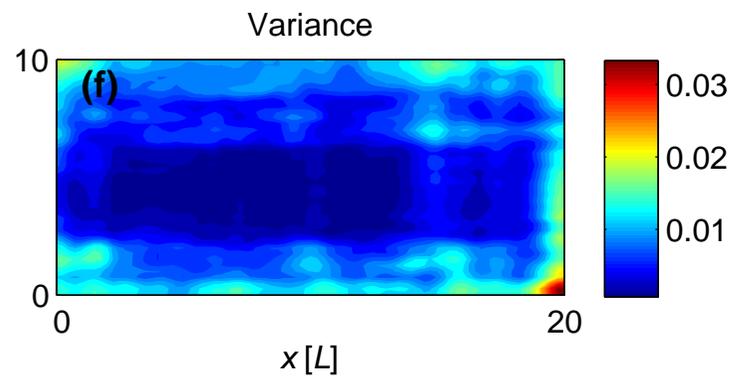